%% file: article-arxiv.tex
\documentclass[11pt]{amsart}
\pdfoutput=1
\usepackage{a4wide}
\usepackage{color}
\usepackage{pgf}
\usepackage[colorlinks, citecolor=green!50!black, linkcolor=red!50!black]{hyperref}
\usepackage{bookmark}
\parskip=\baselineskip

\let\oldlist=\list
\newlength\oldparskip
\def\list#1#2{\oldparskip=\parskip\parskip=0\baselineskip
\oldlist{#1}{#2}\parskip=0\baselineskip}
\let\oldendlist=\enditemize
\def\enditemize{\oldendlist\vspace*{-\oldparskip}}

\renewcommand{\address}[1]{\thanks{#1}}
\renewcommand{\email}[1]{\thanks{email: \texttt{#1}}}
\usepackage{todonotes}

\input{sdpDefinitions}

\begin{document}
\title{Type II$_1$ factors with arbitrary countable endomorphism group}
\stepcounter{footnote}
\author{Steven Deprez}
\thanks{Postdoc at the university of Copenhagen (from September 2011)}
\thanks{Partially supported by ERC Advanced Grant no. OAFPG 247321}
\thanks{Supported by the Danish National Research Foundation (DNRF) through the Centre for Symmetry and Deformation}
\address{Department of mathematics, Copenhagen university, Universitetsparken 5, 2500 O, Copenhagen)}
\email{sdeprez@math.ku.dk}

\begin{abstract}
In \cite{Ioana:vNsuperrigidity}, Ioana introduced three new invariants of type II$_1$
factors: the one-sided fundamental group, the endomorphism semigroup
and the set of right-finite bimodules. In \cite{Ioana:vNsuperrigidity}, he does not
provide many computations of these invariants. In particular, the
question whether these invariants can be trivial is left open. We give
an explicit example of a type II$_1$ factor for which all three
invariants are trivial. More generally, for any countable
left-cancellative semigroup $G$, we construct a type II$_1$ factor $M$
whose endomorphism semigroup is precisely $G$.
\end{abstract}

\maketitle
\section*{Introduction and Overview of the paper}
In \cite[section 10.(II)]{Ioana:vNsuperrigidity}, Ioana introduced three new invariants of type II$_1$
factors, but he provided few concrete computations. Here we provide a
large class of type II$_1$ factors where we can compute these
invariants. The invariants in question are ``one-sided versions'' of
three classical invariants. Let $M$ be a type II$_1$ factor. Then the
one-sided fundamental group $\Fundg_s(M)$ is defined to be
\[\Fundg_s(M)=\{t\in\IRpos\mid \text{ there is a normal injective
}\ast\text{-homomorphism }\varphi:M\rightarrow M^t\}.\]
Observe that this set contains $1$, and is closed under multiplication and under taking
sums. In fact, it is closed under taking infinite sums. This implies
that $\Fundg_s(M)=\IRpos$ whenever $(0,1)\cap\Fundg_s(M)\not=\emptyset$.
In particular, whenever the fundamental group of $M$ is non-trivial,
it follows that the one-sided fundamental group is all of
$\IRpos$. Similarly, it follows that $\Fundg_s(\Lg(\FG_n))=\IRpos$ for
all $n\in\IN$, because $\Lg(\FG_n)\subset\Lg(\FG_{n+1})=\Lg(\FG_n)^t$
where $t=\sqrt{\frac{n-1}n}<1$. Our examples are all on the other side
of the spectrum: they satisfy $\Fundg_s(M)=\IN$.

The second invariant that Ioana introduced is the one-sided version on
the outer automorphism group. This is called the endomorphism
semigroup $\End(M)$. It is the set of all normal injective
$\ast$-homomorphisms $\varphi:M\rightarrow M$, and two such
$\ast$-homomorphisms $\varphi_1,\varphi_2$ are identified if there is
a unitary $u\in M$ such that $\varphi_1=\Ad_u\circ\varphi_2$.
This way, it is clear that $\End(M)$ is a unital semigroup. But it
does not have to be left nor right cancellative. For example,
$\End(R)$, where $R$ is the hyperfinite II$_1$ factor, is neither left
nor right cancellative. This is easy to see explicitly: remember that
$R\cong R\otimes R\cong (R\otimes
R)\rtimes(\IZ/2)$, where $\IZ/2$ acts outerly on $R\otimes R$ by
swapping the components of the tensor product. Write
$\varphi_1,\varphi_2:R\rightarrow R\otimes R$ for the
embeddings that are given by $\varphi_1(x)=x\otimes 1$ and
$\varphi_2(x)=1\otimes x$. Denote by $\psi:R\otimes R\rightarrow
(R\otimes R)\rtimes(\IZ/2)$ the obvious embedding. Then we see that
\begin{align*}
  (\id\otimes\varphi_1)\circ\varphi_1 &=
  (\id\otimes\varphi_2)\circ\varphi_1:R\rightarrow R\otimes R\otimes R
  &\text{ but }
  \id\otimes\varphi_1&\not=\id\otimes\varphi_2\text{ in }\End(R)\\
  \psi\circ\varphi_1&=\psi\circ\varphi_2\text{ in }\End(R)&\text{ but
  }
  \varphi_1&\not=\varphi_2\text{ in }\End(R).
\end{align*}

Even though in general $\End(M)$ does not have to be left
cancellative, our examples will be.
We show that every countable left-cancellative unital semigroup
appears as the endomorphism semigroup of some type II$_1$ factor. In
particular, we find a type II$_1$ factor with trivial endomorphism
semigroup. This solves Ioanas question for an example of such a type
II$_1$ factor.

In fact, we show even more. Ioana introduced a third invariant that
contains both the one-sided fundamental group and the endomorphism
semigroup. This is the set $\RFBimod(M)$ of all $M$-$M$ bimodules $H$
that have finite dimension as a right $M$-module, up to isomorphism of
$M$-$M$ bimodules. This set is closed under the Connes tensor product
and under finite direct sums. It is even closed under infinite direct
sums, provided that the dimensions (i.e.\ the right dimensions over
$M$) form a convergent series.

This invariant contains the previous two invariants. The set
of all right dimensions of $M$-$M$ bimodules is precisely the
one-sided fundamental group. Moreover, the sum and product in $\IRpos$
correspond to the direct sum and the Connes tensor product. The
endomorphism semigroup corresponds precisely to the set of all $M$-$M$
bimodules with right dimension equal to $1$, and the product in
$\End(M)$ corresponds to the Connes tensor product in $\RFBimod(M)$.
We give an example of a type II$_1$ factor for which $\RFBimod(M)$ is
as small as possible, i.e. all right-finite $M$-$M$ bimodules are
trivial bimodules (direct sums of $\Lp^2(M)$).

The results in this paper are based on Popas deformation/rigidity
theory. More precisely, we combine techniques and results from
\cite{Popa:StrongRigidity1,Popa:StrongRigidity2}, \cite{IoanaPetersonPopa:AFP},
\cite{PopaVaes:ActionsOfFinfty,PopaVaes:FundamentalGroupGeneral},
\cite{IoanaPopaVaes:vNsuperrigidity}
and \cite{PopaVaes:CartanFree,PopaVaes:CartanHyperbolic} in order to reduce the computation of
$\End(M)$ to a problem in ergodic theory.

The ergodic-theoretic problem is the following. Let $\Lambda\actson
(Y,\nu)$ be an ergodic probability measure preserving (p.m.p.) action
of a not necessarily countable group. Another
probability measure preserving action $\Lambda\actson (Z,\eta)$ of the
same group is
said to be a factor of $\Lambda\actson (Y,\nu)$ if there is a
p.m.p.\ quotient map $\Delta:Y\rightarrow Z$ such that $\Delta(\lambda
y)=\lambda\Delta(y)$ for almost all $y\in Y$ and this for all
$\lambda\in\Lambda$. The map $\Delta$ is called a factor map. We
denote by $\Factor(\Lambda\actson(Y,\nu))$ the set of all factor maps
from $(Y,\nu)$ to itself. Composition of factor maps defines a
semigroup operation, and the identity map is the identity element for
this operation. This way, we see that $\Factor(\Lambda\actson
(Y,\nu))$ is a right-cancellative semigroup.

Given an ergodic action $\Lambda\actson (Y,\nu)$, we construct a type
II$_1$ factor $M_Y$ such that
$\End(M_Y)=\Factor(\Lambda\actson(Y,\nu))^{\op}$. It is easy to see
that every compact group $G$ is $G=\Factor(G\actson (G,h))$ where $h$
denotes the Haar measure on $G$. Hence $G$ appears also as the
endomorphism semigroup of some type II$_1$ factor $M$. In these cases,
all endomorphisms of $M$ are in fact isomorphisms. Every compact
right-cancellative unital semigroup is automatically a group, so this
observation covers whole compact right-cancellative case.

The discrete case is more interesting. Here we have examples of
semigroups that are not groups. Already the semigroup of natural
numbers with addition form such a semigroup. In section
\ref{sect:semigroups}, we show that every countable right-cancellative
semigroup appears as the factor semigroup of an ergodic p.m.p.\ action
$\Lambda\actson (Y,\nu)$. Hence every countable left-cancellative
semigroup appears as $\End(M)$ for some type II$_1$ factor $M$.

Given $\Lambda\actson (Y,\nu)$, we construct $M_Y$ as follows. We can
consider $\Lambda$ (in fact, a quotient of $\Lambda$) as a subgroup of
$\Autmp(Y,\nu)$. Observe that $\Factor(\Lambda\actson(Y,\nu))$ only
depends on the closure of $\Lambda$ in the usual Polish topology on
$\Autmp(Y,\nu)$. Hence we can replace $\Lambda$ by a countable dense
subgroup without changing $\Factor(\Lambda\actson(Y,\nu))$. From now
on we assume that $\Lambda$ is countable. Let $\Gamma_1$ be a
hyperbolic property (T) group with trivial endomorphism semigroup. Let
$\Sigma\subset\Gamma_1$ be an
amenable subgroup. Consider $\Gamma=\Gamma_1\free_\Sigma
(\Sigma\times\Lambda)$. Let $\Gamma\actson I$ be an action of $\Gamma$
on a countable set.

Let $(X_0,\mu_0)$ be an atomic probability space with unequal weights, and
consider the generalized Bernoulli action $\Gamma\actson
(X,\mu)=(X_0,\mu_0)^I$. Consider the obvious quotient map
$\pi:\Gamma\rightarrow\Lambda$. Following
\cite{PopaVaes:ActionsOfFinfty,PopaVaes:FundamentalGroupGeneral}, we
define an action of $\Gamma$ on $X\times Y$ by the formula
$g(x,y)=(gx,\pi(g)y)$. Then we set $M_Y=\Lp^\infty(X\times
Y)\rtimes\Gamma$.

Let us give an idea why $\End(M_Y)=\Factor(\Lambda\actson
(Y,\nu))^{\op}$. One inclusion is easy. Given $\Delta\in\Factor(\Lambda\actson (Y,\nu))$, we
define an embedding $\varphi_\Delta:M_Y\rightarrow M_Y$ by the formula
$\varphi((a\otimes b)u_g)=(a\otimes\Delta_\ast(b))u_g$, where
$\Delta_\ast(b)=b\circ\Delta$ for every function
$b\in\Lp^\infty(Y)$. The application $\Delta\mapsto\varphi_\Delta$
embeds $\Factor(\Lambda\actson Y)^{\op}$ into $\End(M_Y)$.

Now, let $\varphi:M_Y\rightarrow M_Y$ be an endomorphism of
$M_Y$. We want to show that $\varphi=\varphi_\Delta$ up unitary
conjugacy. Denote $A=\Lp^\infty(X)$ and $B=\Lp^\infty(Y)$. Techiques
from \cite{Popa:StrongRigidity1,Popa:StrongRigidity2} show that
$\varphi(B\rtimes\Gamma)\subset B\rtimes\Gamma$, up to a
unitary. This result depends crucially on the fact the $\Gamma_1$ has
property (T) while the action $\Gamma\actson X$ is a generalized
Bernoulli action. Similarly,
techniques from \cite{IoanaPetersonPopa:AFP} show that $\varphi((A\otimes
B)\rtimes\Gamma_1)\subset (A\otimes B)\rtimes\Gamma_1$ up to unitary
conjugacy. In fact, we can assume that both unitaries are the
same. This result uses the facts that $\Gamma_1$ has property (T)
while $\Gamma$ is an amalgamated free product. Because $\Gamma_1$ is
hyperbolic, \cite{PopaVaes:CartanFree,PopaVaes:CartanHyperbolic} shows
that $\varphi(A)$ can not be in
$B\rtimes\Gamma_1=B\otimes\Lg(\Gamma_1)$. Then \cite[theorem
  5.1]{IoanaPopaVaes:vNsuperrigidity} shows that $C=\varphi(A)^\prime\cap (A\otimes
B)\rtimes\Gamma_1$ embeds into $A\otimes B$. Now we apply
\cite[theorem 6.1]{IoanaPopaVaes:vNsuperrigidity} to conclude that $\varphi(A\otimes B)\subset
A\otimes B$. Moreover, the endomorphism
$\varphi:B\otimes\Lg(\Gamma_1)\rightarrow B\otimes\Lg(\Gamma_1)$ is
described in the following way. We can consider every element in
$B\otimes\Lg(\Gamma_1)$ as a map from $Y$ to $\Lg(\Gamma_1)$. There is
a field of group endomorphisms $\delta_y:\Gamma_1\rightarrow\Gamma_1$
($y\in Y$) such that $\varphi(u_g)(y)=u_{\delta_y(g)}$. All of these
are inner, so we can assume that they are all trivial. Now we know
that $\varphi(A\otimes B)\subset A\otimes B$ and $\varphi(u_g)=u_g$
for all $g\in\Gamma_1$. For a good choice for the action
$\Gamma\actson I$, a direct computation shows that in fact
$\varphi=\varphi_{\Delta}$ for some factor map
$\Delta\in\Factor(\Lambda\actson Y)$.

Of course, in the above idea of the proof, we have been ignoring a lot
of technical conditions. A more precise statement and proof are given
in \ref{sect:main}. In section \ref{sect:prelim}, we remind the reader
of some well-known results that are crucial for this paper. In section
2, we extend \cite[theorem 5.1 and 6.1]{IoanaPopaVaes:vNsuperrigidity}
to our setting. Section 3 introduces two properties of groups that are
crucial in the next section. There we show our main result, theorem
\ref{thm:main}. In order to apply that main theorem, we need to give
an example of a group $\Gamma$ and an action $\Gamma\actson I$ that
satisfies the conditions of theorem \ref{thm:main}. This is not very
hard, but it is technical. Section \ref{sect:ex} is devoted to such an
example. Finally, in section \ref{sect:semigroups}, we show that all
countable right-cancellative semigroups appear as $\End(M)$ for some
type II$_1$ factor $M$.

\section{Preliminaries and Notations}
\label{sect:prelim}
\subsection{Relatively weakly mixing actions}
Relative weak mixing plays a crucial role in the proof of
theorem \ref{thm:inter-crit}. This property was introduced by Furstenberg in
\cite{Furstenberg:DiagonalMeasures} and Zimmer in
\cite{Zimmer:GenDiscreteSpectrum, Zimmer:Extensions}, in the case of
actions on probability spaces. In \cite{Popa:CocycleSuperrigidityMalleable}, Popa generalized
this to actions on von Neumann algebras.

\begin{definition}[see {\cite[lemma 2.10]{Popa:CocycleSuperrigidityMalleable}}]
  Let $D\subset (B,\tau)$ be an inclusion of finite von Neumann
  algebras. Assume that a countable group $\Gamma$ acts trace-preservingly on
  $B$ and leaves $D$ globally invariant. Denote the action by
  $\alpha$. We say that $\Gamma$ acts
  weakly mixingly on $B$ relative to $D$ if one of the following
  equivalent conditions holds.
  \begin{enumerate}
  \item There exists a sequence of group elements $(g_n)_n$ in
    $\Gamma$ such that
    \[\norm{\E_D(x\alpha_{g_n}(y))}_2\rightarrow 0\qquad\text{for all
    }x\in B,y\in B\ominus D.\]
  \item Every $\Gamma$-invariant positive element $a$ with finite
    trace in the basic construction $\langle B,e_D\rangle$, must be $a\in e_D D e_D$.
  \item The action of $\Gamma$ on $\Lp^2(B)\otimes_D\Lp^2(B)$ is
    ergodic relative to $\Lp^2(D)$. I.e. all $\Gamma$-invariant
    vectors $\xi\in\Lp^2(B)\otimes_D\Lp^2(B)$ are $\xi\in \Lp^2(D)$.
  \item The only right $D$-submodules of $\Lp^2(B)$ that have finite
    dimension over $D$ and are $\Gamma$-invariant, are contained in
    $\Lp^2(B)$.
  \end{enumerate}
\end{definition}
\begin{proof}[proof of equivalence of these conditions]
  For a proof that conditions (1) and (2) are equivalent, we refer to
  \cite[lemma 2.10]{Popa:CocycleSuperrigidityMalleable}. Equivalence of conditions (2)
  and (3) follows from the fact the
  $\Lp^2(\langle B,e_D\rangle)=\Lp^2(B)\otimes_D\Lp^2(B)$. Condition (2) and (4)
  are equivalent because the $D$-submodules of $\Lp^2(B)$ correspond
  1-to-1 with the projections in $\langle B,e_D\rangle$. The dimension
  of the submodule is precisely the trace of the corresponding
  projection, and the submodule is $\Gamma$-invariant if and only if
  its corresponding projection is.
\end{proof}

In the case where $B$ is abelian, we recover the classical definition
of relative weak mixing:
\begin{observation}
  Let $\Gamma\actson (X,\mu)$ be a p.m.p. action and suppose that
  $p:(X,\mu)\rightarrow (Y,\nu)$ is a factor of this action. Then
  $\Lp^\infty(Y,\nu)\subset\Lp^\infty(X,\mu)$ is a $\Gamma$-invariant
  von Neumann subalgebra. Then the following are equivalent:
  \begin{enumerate}
  \item $\Gamma\actson \Lp^\infty(X,\mu)$ is weakly mixing relative to
    $\Lp^\infty(Y,\nu)$, in the sense defined above.
  \item $\Gamma\actson (X,\mu)$ is weakly mixing relative to
    $p:X\rightarrow Y$, in the classical sense. I.e. the diagonal
    action of $\Gamma$ on $X\times_Y X$ is ergodic.
  \end{enumerate}
\end{observation}

\subsection{Generalized co-induced actions}
\label{subsect:prelim:gencoind}
We introduced generalized co-induced actions in \cite{Deprez:Fundg}. Here we
repeat the construction and generalize some properties.

\begin{definition}
  \label{def:gen-co-ind}
  Let $\Lambda\actson I$ be an action of a countable group on a
  countable set. Let $\omega:\Lambda\times I\rightarrow \Lambda_0$ be
  a cocycle. Suppose that $\Lambda_0$ acts probability measure
  preservingly on $(Y_0,\nu_0)$. Define an action of $\Lambda$ on
  $(Y,\nu)=(Y_0,\nu_0)^I$ by the formula
  $(gy)_i=\omega(g,g^{-1}i)y_{g^{-1}i}$.
  This action is called the generalized co-induced action of
  $\Lambda_0\actson (Y_0,\nu_0)$, with respect to $\omega$.
\end{definition}
\begin{lemma}
  \label{lem:gen-co-ind}
  Let $\Lambda\actson I$ be an action of a countable group on a
  countable set, and let\\$\omega:\Lambda\times I\rightarrow\Lambda_0$
  be a cocycle. Suppose that $\Lambda_0\actson (Y_0,\nu_0)$ is a
  probability measure preserving action. Consider the generalized
  co-induced action $\Lambda\actson (Y,\nu)=(Y_0,\nu_0)^I$.
  \begin{itemize}
  \item If all orbits of $\Lambda\actson I$ are infinite, then the
    generalized co-induced action $\Lambda\actson(Y,\nu)$ is weakly mixing.
  \item Suppose that $\Lambda\actson I$ and $\omega$ satisfy the
    following three conditions.
    \begin{itemize}
    \item $\Lambda$ acts transitively on $I$.
    \item There exists an $i\in I$ such that (or equivalently, for all
      $i\in I$) $\omega$ maps the set
      $\Stab\{i\}\times\{i\}$ surjectively onto $\Lambda_0$
    \item There exists an $i\in I$ such that (or equivalently, for all
      $i\in I$) the subgroup $S_i=\{g\in\Lambda\mid gi=i\text{ and
      }\omega(g,i)=e\}$ acts with infinite orbits on $I\setminus\{i\}$.
    \end{itemize}
    Then every measurable $\Lambda$-invariant map $f:Y\rightarrow Y$
    is of the form $f(x)_i=f_0(x_i)$, where $f_0:Y_0\rightarrow Y_0$
    is a measurable $\Lambda_0$-invariant map.
  \end{itemize}
\end{lemma}
\begin{proof}
  The first point follows in the same way as for generalized Bernoulli
  actions. See for example \cite[proposition
  2.3]{PopaVaes:StrongRigidity}.

  For the second point, it is clear that every map of the given form
  is measurable and $\Lambda$-invariant. On the other hand, let
  $f:Y\rightarrow Y$ be a $\Lambda$-invariant map. Fix $i\in I$ and
  consider the composition $f_i:Y\rightarrow Y_0$ of $f$ with the
  quotient onto the $i$-th component of $Y$. Observe that $f_i$ is
  $S_i$-invariant, because $S_i$ does not act on the $i$-th
  component. But by the first point, $S_i$ acts ergodically on
  $Y_0^{I\setminusb\{i\}}$. So we see that $f(x)_i=f_i(x)=f_0(x_i)$
  for some measurable map $f_0:Y_0\rightarrow Y_0$, but only for the
  one $i\in I$ we fixed.

  Moreover, observe that $f_0(\omega(g,i)x_0)=\omega(g,i)f_0(x_0)$ for
  all $g\in\Stab\{i\}$. Since $\omega$ maps $\Stab{i}\times\{i\}$ onto
  $\Lambda_0$, we find that $f_0$ is $\Lambda_0$-equivariant. Let
  $j\in I$ be another element. We find $g\in\Lambda$ such that
  $gj=i$. Then we compute that
  \[f(x)_j=\omega(g^{-1},j)(gf(x))_i=\omega(g,i)^{-1}f_0((gx)_i)=f_0(x_j).\]
\end{proof}

\section{Generalizing two results from \cite{IoanaPopaVaes:vNsuperrigidity}}
Our main theorem depends on the results from sections 5 and 6 from
\cite{IoanaPopaVaes:vNsuperrigidity}. But in fact, we need a slightly more general statement of
these two results. The proof is mainly an application of the
direct integral decomposition of a von Neumann algebra. Nevertheless,
we think it is worthwhile to give a careful statement and a short
proof.

\begin{theorem}[{a version of \cite[section
    5]{IoanaPopaVaes:vNsuperrigidity}}]
  \label{thm:inter-cartan}
  Let $\Gamma$ act on a countable set $I$ in such a way that there is
  a number $\kappa\in\IN$ such that $\Stab\cF$ is finite whenever
  $\nelt{\cF}\geq \kappa$. Choose a standard probability space
  $(X_0,\mu_0)$. Suppose that $(B,\tau)$ is a finite type I von Neumann
  algebra. Write $A=\Lp^\infty(X_0^I)$ and consider the von Neumann
  algebra $M=(A\rtimes\Gamma)\otimes B$.

  Let $p\in\Lg(\Gamma)\otimes B$ be a projection.
  Let $D\subset pMp$ be an abelian subalgebra. Write $\cG$ for the
  normalizer of $D$ inside $pMp$. Denote the intersection
  $\cG_0=\cG\cap\Unitary(p(\Lg(\Gamma)\otimes B)p)$. Assume that
  \begin{itemize}
  \item $D$ does not embed into $B$ inside $M$.
  \item $\cG^{\prime\prime}$ does not embed into
    $(A\rtimes\Stab\{i\})\otimes B$ inside $M$, for any $i\in I$.
  \item $\cG^{\prime\prime}$ does not embed into
    $\Lg(\Gamma)\otimes B$ inside $M$.
  \item $\cG_0^{\prime\prime}$ does not embed into
    $(\Lg\Stab\{i\})\otimes B$ inside $\Lg(\Gamma)\otimes B$
    for any $i\in I$.
  \end{itemize}
  Then we get that $C=D^\prime\cap pMp\fembeds_M A\otimes B$.
\end{theorem}
\begin{proof}
  Set $\widetilde D=\Centre(C)$, and observe that $\widetilde D$ is
  still an abelian subalgebra of $pMp$ that is still normalized by
  $\cG$, and $D\subset \widetilde D$. Hence $\widetilde D$ still satisfies the four conditions
  above. But we also get that $p\Centre(B)\subset \widetilde D$.

  We write $\Centre(B)=\Lp^\infty(Y,\nu)$ for some standard measure
  space $(Y,\nu)$. Then we can take the direct integral decomposition
  $B=\int^\oplus_Y B_y d\nu(y)$. Likewise we can decompose
  $p=\int^\oplus_Y p_y d\nu(y)$, where each $p_y\in
  \Lg(\Gamma)\otimes B_y$.

  We decompose $\widetilde D$ and $C$ into $\widetilde D=\int^\oplus_Y D_y d\nu(y)$ and
  $C=\int^\oplus_Y C_y d\nu(y)$ respectively. Denote
  $M_y=(A\rtimes\Gamma)\otimes B_y$, and observe that $C_y$ is the
  relative commutant of $D_y$ inside $p_yM_yp_y$. Each unitary $u$ in
  $\cG$ decomposes intro a direct integral of unitaries $u_y$, each of
  which normalizes $D_y$.

  All in all, we see that the inclusion $D_y\subset p_yM_yp_y$
  satisfies the conditions of \cite[theorem 5.1]{IoanaPopaVaes:vNsuperrigidity}. We obtain that
  $C_y\fembeds_{M_y}A\otimes B_y$ for almost all $y\in Y$. Hence it
  follows that $C\fembeds_M A\otimes B$.
\end{proof}

We also want to give a similar variant to \cite[theorem 6.1]{IoanaPopaVaes:vNsuperrigidity} But
for our main theorem, we need a slightly more general version: using
the notations from \cite[section 6]{IoanaPopaVaes:vNsuperrigidity}, we can not assume that
$(\Ad_{\gamma(s)})_{s\in\Lambda}$ acts weakly mixingly on
$\Centre(C)$.  Instead, we can only assume that the action is weakly mixing
relative to a discrete subalgebra $D\subset \Centre(C)$. This is not a
hard generalization, but we have to adapt the statement of the theorem
slightly. In order to simplify the statement of theorem \ref{thm:inter-crit}, we
incorporate \cite[corollary 6.2.1]{IoanaPopaVaes:vNsuperrigidity}.
\begin{theorem}[{our variant of \cite[theorem 6.1 and corollary 6.2.1]{IoanaPopaVaes:vNsuperrigidity}}]
  \label{thm:inter-crit}
  Let $\Gamma\actson (X,\mu)$ be a free, ergodic and p.m.p.\
  action. Let $(B,\tau)$ be a finite type I von neumann algebra. Write
  $A=\Lp^\infty(X)$ and consider the von Neumann algebra
  $M=(A\rtimes\Gamma)\otimes B$. Let $p\in \Lg(\Gamma)\otimes B$ be a
  projection with finite trace.

  Let $C\subset pMp$ be a von Neumann subalgebra and suppose that
  $\gamma:\Lambda\rightarrow\Unitary(p(\Lg(\Gamma)\otimes
  B)p)\cap\vnNorm_{pMp}(C)$ is a group morphism satisfying the
  following conditions. 
  \begin{itemize}
  \item $\Lambda$ does not have any non-trivial finite dimensional
    unitary representations.
  \item $\gamma(\Lambda)^{\prime\prime}$ does not embed into any
    $\Lg(\Centr\{g\})\otimes B$ for any $e\not=g\in\Gamma$.
  \item Consider the action of $\Lambda$ on $\Centre(C)$ by conjugating with
    $\gamma(\Lambda)$. We assume that this action is weakly mixing
    relative to $D=\Centre(C)\cap p(\Lg(\Gamma)\otimes B)p$.
  \item $\Centre(C)^\prime\cap pMp=C$ and $C\fembeds A\otimes B$.
  \end{itemize}

  Then we have the following
  \begin{itemize}
  \item a partial isometry $v\in \Lg(\Gamma)\otimes B\otimes
    \Bounded(\IC,\ell^2(\IN))\otimes\Bounded(\IC,\ell^2(\IN))$ with
    left support equal to $p$ and with right
    support $q=v^\ast v$ inside $\widetilde
    B=B\otimes\Bounded(\ell^2(\IN))\otimes\ell^\infty(\IN)$, and
  \item a group morphism $\delta:\Lambda\rightarrow \cG$ where
    $\cG\subset\Unitary(\Lg(\Gamma)\otimes \widetilde B)$
    is the group
    \[\cG=\left\{\left.\sum_{g\in\Gamma}p_gu_g\,\right\vert\, p_g\in\Centre(q\widetilde
      Bq)\text{ are projections with }\sum_gp_g=q\right\},\]
  \end{itemize}
  such that
  \[v^\ast Cv=q(A\otimes\widetilde B)q\text{ and
  }v^\ast\gamma(s)v=\delta(s)\text{ for all }s\in\Lambda.\]
\end{theorem}
\begin{proof}
  Write $\Centre(B)=\Lp^\infty(Y,\nu)$ for some measure space
  $(Y,\nu)$. Decompose $B$ as the direct integral of $(B_y)_{y\in Y}$.
  Observe that $p\Centre(B)\subset D\subset \Centre(C)$. Hence we can
  decompose both $D$ and $C$ as a direct integral over $Y$ of von
  Neumann algebras $D_y$ and $C_y$. Write $M_y=(A\rtimes\Gamma)\otimes
  B_y$ and observe that $M$ is the direct integral of the $M_y$. We
  decompose $p$ as the direct integral of the projections $p_y$.

  Observe that $C_y\fembeds_{M_y} A\otimes B_y$ and that
  $\Centre(C_y)^\prime\cap p_yM_yp_y=C_y$ almost everywhere. Moreover, we can
  consider the $\gamma(s)$ as measurable maps $y\mapsto\gamma(s,y)$ from $Y$ into the group
  of unitary elements of $p_y(\Lg(\Gamma)\otimes B_y)p_y$. We remark that the
  $\gamma(s,y)$ normalize $C_y$. Consider the action of $\Lambda$ on
  $\Centre(C_y)$ that is given by conjugation with $\gamma(s,y)$. This
  action is still weakly mixing relative to $D_y$.

  But $D_y$ embeds fully into $A\otimes B_y$, while $D_y$ is contained
  in $\Lg(\Gamma)\otimes B_y$. Hence it fully embeds into $B_y$ (see
  lemma \ref{lem:gmc-group} below). Since $B_y$ is a finite type I factor, it follows
  that $D_y$ is an abelian discrete von Neumann algebra. Since
  $\Lambda$ does not have non-trivial finite dimensional unitary
  representations, it can not act trace-preservingly on such a von
  Neumann algebra, unless the action is trivial. So we find a
  countable number of projections $p_{y,n}\in D_y$ such that $\sum_n
  p_{n,y}=p_y$ and $D_yp_{n,y}=\IC p_{n,y}$. In other words, $\Lambda$
  acts weakly mixingly on $\Centre(C_y)p_{n,y}$.

  We can now apply \cite[theorem 6.1 and corollary 6.2.1]{IoanaPopaVaes:vNsuperrigidity} to the
  inclusion $C_yp_{n,y}\subset p_{n,y}M_yp_{n,y}$. We find a partial
  isometry $v_{n,y}\in \Lg(\Gamma)\otimes B_y\otimes\Bounded(\IC,\ell^2(\IN))$ with left
  support equal to $p_{n,y}$ and with right support
  $q_{n,y}=v_{n,y}^\ast v_{n,y}\in B_y\otimes\Bounded(\ell^2(\IN))$,
  and such that
  \[v_{n,y}^\ast C_y v_{n,y}=q_{n,y}(A\otimes
  B_y\otimes\Bounded(\ell^2(\IN)))q_{n,y}.\]
  Moreover, we find a group morphism
  $\delta_{y,n}:\Lambda\rightarrow\Gamma$ and a finite-dimensional
  unitary representation
  $\pi_{n,y}:\Lambda\rightarrow\Unitary(p_{n,y}(B_y\otimes\Bounded(\ell^2(\IN)))p_{n,y})$
  such that
  \[v_{n,y}^\ast \gamma(s,y)v_{n,y}=\pi_{n,y}(s)\otimes
  u_{\delta(s)}\text{ for all }s\in\Lambda.\]
  But we assumed that $\Lambda$ does not have any non-trivial finite
  dimensional unitary representations, so we see that
  $\pi_{n,y}(s)=p_{n,y}$.

  In fact, reading the proof of \cite[theorem 6.1]{IoanaPopaVaes:vNsuperrigidity} carefully, we
  see that we can do all this in such a way that the $v_{n,y}$ depend
  measurably on $y$. Hence we can consider the partial isometry
  \[v=\sum_n\left(\int_Y^\oplus v_{n,y}d\nu(y)\right)\otimes e_{1,n}\in
  \widetilde N\Lg(\Gamma)\otimes B\otimes\Bounded(\IC,\ell^2(\IN))\otimes\Bounded(\IC,\ell^2(\IN)).\]
  This partial isometry has left support $p$ and its right support is
  given by $q=\int_Y^\oplus \sum_n q_{n,y}\otimes
  e_{n,n}d\nu(y)\in\widetilde B$. A direct computation shows that
  \[v^\ast Cv=q(A\otimes \widetilde B)q.\]

  The measureable field $\delta(s,y,n)=\delta_{n,y}(s)$ of group
  morphisms satisfies the condition
  \[v^\ast \gamma(s) v=q\sum_g
  \chara_{\{(y,n)\mid\delta(s,y,n)=g\}}\otimes u_{g},\]
  so we see that $v^\ast \gamma(s)v\in \cG$ for all $s\in\Lambda$.
\end{proof}

\begin{lemma}
  \label{lem:gmc-group}
  Let $A,B$ be finite von Neumann algebras with
  trace-preserving actions of a countable group $\Gamma$. Consider
  $M=(A\otimes B)\rtimes\Gamma$. Let $D\subset p(B\rtimes\Gamma)p$ be
  a von Neumann subalgebra of some corner of $B\rtimes\Gamma$. If $D$
  embeds into $A\otimes B$ inside $M$, then $D$ already embeds into
  $B$ inside $B\rtimes\Gamma$.
\end{lemma}
\begin{proof}
  Suppose that $D$ did not embed into $B$ inside
  $B\rtimes\Gamma$. By definition, we find a sequence $(v_n)_n$ of
  unitaries in $D$ such that
  \[\norm{\E_B(xv_ny)}_2\rightarrow 0\text{ for all }x,y\in
  B\rtimes\Gamma.\]
  We want to show that $D$ does not embed into $A\otimes B$, inside
  $M$. So we want to show that
  \[\norm{\E_{A\otimes B}(xv_ny)}_2\rightarrow 0\text{ for all }x,y\in
  M.\]
  By Kaplansky's density theorem, we can assume that $x=(a_1\otimes
  b_1)u_{g}$ and $y=u_h(a_2\otimes b_2)$ for some $a_1,a_2\in A$ and
  $b_1,b_2\in B$ and $g,h\in\Gamma$. We compute that
  \begin{align*}
    \norm{\E_{A\otimes B}(xv_ny)}_2
    &=\norm{(a_1\otimes b_1)\E_{A\otimes B}(u_gv_nu_h)(a_2\otimes
      b_2)}_2\\
    &\leq\norm{a_1}\norm{a_2}\norm{b_1}\norm{b_2}\norm{\E_B(u_gv_nu_h)}_2\rightarrow
    0,
  \end{align*}
  as required.
\end{proof}

\section{Anti-(T) groups and groups with small normalizers}
It is well-known that a group $\Gamma$ that has the Haagerup property
does not contain an infinite property (T) subgroup. Slighly more
generally, for any p.m.p.\ action $\Gamma\actson(X,\mu)$, we know that
the corresponding group-measure space von Neumann algebra
$\Lp^\infty(X)\rtimes\Gamma$ does not contain a diffuse von Neumann
subalgebra with property (T), see \cite{Popa:Betti}.
\begin{definition}[anti-(T) group]
  \label{def:anti-T}
  We say that a group $\Gamma$ is anti-(T) if for every trace
  preserving action $\Gamma\actson (A,\tau)$ on an amenable von
  Neumann algebra $A$ and for every projection $p\in A\rtimes\Gamma$,
  the von Neumann algebra $p(A\rtimes\Gamma)p$ does not contain a
  diffuse von Neumann subalgebra with property (T).
\end{definition}
This definition differs from the notion of an anti-(T) group
that was introduced in
\cite{HoudayerPopaVaes:AppVNsuperrigidity}. Every group that has the
anti-(T) property in the
\cite{HoudayerPopaVaes:AppVNsuperrigidity}-sense is anti-(T) in our
sense, but not the other way around. The advantage of our notion is
that it is stable under arbitrary amalgamated free products.
\begin{lemma}
  If $\Gamma_1,\Gamma_2$ are anti-(T) groups, then
  $\Gamma=\Gamma_1\free_\Sigma\Gamma_2$ is still anti-(T).
\end{lemma}
\begin{proof}
  Let $\Gamma\actson (A,\tau)$ be a trace-preserving action on an
  amenable von Neumann algebra. Let $p\in A\rtimes\Gamma$ be a
  projection and suppose that $Q\subset p(A\rtimes\Gamma)p$ is a
  property (T) subalgebra. By \cite[theorem
  5.1]{IoanaPetersonPopa:AFP}, we know that $Q\embeds
  A\rtimes\Gamma_i$ for $i=1$ or $2$. In particular, there is a
  $\ast$-homomorphism $\theta:Q\rtimes q(A\rtimes\Gamma_i)^nq$ for
  some $n\in \IN$ and $q\in (A\rtimes\Gamma_i)^n$. Since $\Gamma_i$
  was anti-(T), it follows that $Q$ is not diffuse.
\end{proof}

In \cite{OzawaPopa:UniqueCartan}, Ozawa and Popa show that the free
groups have the following property: the normalizer of every diffuse
abelian subalgebra $B\subset \Lg(\FG_n)$ is amenable. This
result was generalized later in \cite{PopaVaes:CartanFree}. There Popa
and Vaes
show that for every
trace-preserving action $\FG_n\actson (A,\tau)$ on a finite amenable von
Neumann algebra and every diffuse abelian
subalgebra $B\subset p(A\rtimes\FG_n)p$ of a corner of the crosses
product, we have the following dichotomy. Either $B$ embeds into
$A$ or the normalizer of $B$ is amenable. In
\cite{PopaVaes:CartanHyperbolic}, they show that this property holds
for all hyperbolic groups.

In this section, we introduce a similar but weaker property of
groups. The advantage of our property is that it is implied by the
Haagerup property and that it is stable under taking amalgamated free
products.

We first introduce what we mean when we say that an abelian subalgebra
$B\subset M$ has a large normalizer. In the following, we denote by
$\MatD_n(\IC)$ the algebra of diagonal $n\times n$ matrices with
complex entries.
\begin{definition}
  Let $B\subset M$ be an abelian subalgebra of a finite von
  Neumann algebra $M$. We say that $B$ has a large
  normalizer if $\vnNorm_M(B)$ contains a group $\cG$ that generates a
  diffuse property (T) subalgebra of $M$.
\end{definition}

\begin{definition}[groups with small normalizers]
  We say that a group $\Gamma$ has small normalizers if, for
  every trace-preserving action $\Gamma\actson (A,\Tr)$ on a finite
  amenable von Neumann algebra, and every projection $p\in
  M=A\rtimes\Gamma$, every diffuse abelian
  subalgebra $B\subset pMp$ with large normalizer, embeds into $A$
  inside $M$.
\end{definition}

It follows immediately from \cite[theorem
3.1 and lemma 2.4]{PopaVaes:CartanHyperbolic} that the following
groups have small normalizers.
\begin{itemize}
\item hyperbolic groups
\item lattices in rank one simple Lie groups with finite center
\item Sela's limit groups
\end{itemize}

Any anti-(T) group $\Gamma$ has small normalizers,
simply because no amplification of the crossed product
$A\rtimes\Gamma$ can have a von Neumann subalgebra with property (T).

\begin{theorem}
\label{thm:smallnorm-afp}
The amalgamated free product of groups with small normalizers over an
amalgam that has the anti-(T) property still has small normalizers.
\end{theorem}
\begin{proof}
  Let $\Gamma_1,\Gamma_2$ be two groups with small normalizers, and
  let $\Sigma$ be a common subgroup with the Haagerup
  property. Consider the amalgamated free
  product $\Gamma=\Gamma_1\free_\Sigma\Gamma_2$. Let $\Gamma$ act
  trace preservingly on an amenable finite von Neumann algebra
  $(A,\tau)$. Denote $M=A\rtimes\Gamma$ and take a projection $p\in
  M$. Let $B\subset pMp$ be an abelian subalgebra with a large
  normalizer. Consider a subgroup $\cG\subset \vnNorm_{pMp}(B)$ that
  generates a property (T) subalgebra in $pMp$. Denote by
  $N=B\vee\cG\subset pMp$ the von Neumann subalgebra that is generated
  by $B$ and $\cG$.

  Now we use techniques from \cite{IoanaPetersonPopa:AFP}. But we use the
  versions as stated in \cite{PopaVaes:vNsuperrigidity}, because these versions are
  more convenient. Consider the word-length deformation as
  defined in \cite[section 2.3]{IoanaPetersonPopa:AFP} (see also \cite[section 5.1]{PopaVaes:vNsuperrigidity}): we
  define completely positive maps $m_\rho:M\rightarrow M$ by
  $m_\rho(au_g)=\rho^{\len{g}}au_g$ where $\len{g}$ denotes the
  word-length of $g$ and $\rho$ is a real number between 0 and
  1. These completely positive maps converge to the identity pointwise, as
  $\rho\rightarrow 1$.

  Because $\cG^{\prime\prime}$ has property (T), we know that $m_\rho$
  converges to the identity uniformly in $\norm{\cdot}_2$ on $\cG$. 
  Observe that $\cG^{\prime\prime}$ does not embed into
  $A\rtimes\Sigma$ because $\Sigma$ is anti-(T). So \cite[lemma 5.7]{PopaVaes:vNsuperrigidity} yields a real
  number $0<\rho_0<1$ and a $\delta>0$ such that
  \[\tau(w^\ast m_{\rho_0}(w))>\delta\text{ for all unitaries }w\in B.\]

  By property (T) we find $\rho\geq\rho_0$ such that
  $\norm{v-m_{\rho}(v)}_2<\frac18\delta^2$ for all $v\in\cG$.
  Because $\tau(w^\ast m_{\rho_0}(w))$ increases with $\rho_0$, we can
  assume that $\rho=\rho_0$. Hence we can conclude that\\$\tau(v^\ast
  w^\ast m_{\rho_0}(wv))>\frac12\delta$ for all $v\in \cG$ and
  $w\in\Unitary(B)$.

  Theorem \cite[theorem 4.3]{IoanaPetersonPopa:AFP} shows that $N$ embeds into either
  $A\rtimes\Gamma_1$ or $A\rtimes\Gamma_2$. Without loss of generality,
  we can assume that $N$ embeds into $A\rtimes\Gamma_1$. Hence we
  find a non-zero partial isometry $v\in \MatM_{1,n}(\IC)\otimes M$
  and a $\ast$-homomorphism $\theta:N\rightarrow q(\MatM_n(\IC)\otimes
  A\rtimes\Gamma_1)q$ such that $xv=v\theta(x)$ for all $x\in N$. We
  can assume that $q$ is the support projection of
  $\E_{\MatM_n(\IC)\otimes A}(v^\ast v)$.
  The subalgebra $\theta(B)\subset q((\MatM_n(\IC)\otimes
  A)\rtimes\Gamma_1)q$ still has large normalizer. Because $\Gamma_1$ was assumed to have small
  normalizers, we see that $\theta(B)$ embeds into
  $\MatM_n(\IC)\otimes A$ inside $(\MatM_n(\IC)\otimes
  A)\rtimes\Gamma_1$. It follows that $B$
  embeds into $A$.
\end{proof}

\section{Proof of the main result}
\label{sect:main}
We want to show that all semigroups of the form
$\Factor(\Lambda\actson Y)^{\op}$ appear is $\End(M)$ for some type
II$_1$ factor $M$. We can always replace $\Lambda$ by a countable
dense subgroup of $\Lambda\subset \Aut(Y,\nu)$ without changing
$\Factor(\Lambda\actson Y)$. From now on we assume that $\Lambda$ is
countable. Lemma \ref{lem:no-fd} in section \ref{sect:semigroups} below shows
that we can also assume that $\Lambda$ is anti-(T) and that all
cocycles $\omega:\Lambda\times Y\rightarrow K$ with values in a
compact group are trivial.

Recall the construction of $M$ from the introduction. Let
$\Gamma=\Gamma_1\free_\Sigma(\Sigma\times\Lambda)$ be an amalgamated
free product group. Take an action $\Gamma\actson I$ of $\Gamma$ on a
countable set $I$. Choose a purely atomic base space $(X_0,\mu_0)$
with unequal weights. Set $(X,\mu)=(X_0,\mu_0)^I$ and consider the
generalized Bernoulli action $\Gamma\actson X$. Consider that
canonical quotient morphism $\pi:\Gamma\rightarrow\Lambda$ and define
a new action $\Gamma\actson X\times Y$ by the formula
$g(x,y)=(gx,\pi(g)y)$ for all $g\in\Gamma$ and almost all $(x,y)\in
X\times Y$.

Define $M=\Lp^\infty(X\times Y)\rtimes\Gamma$. Observe that every
$\Delta\in\Factor(\Lambda\actson Y)$ defines an endomorphism
$\varphi_\Delta:M\rightarrow M$ by the formula
$\varphi_\Delta((a\otimes b)u_g)=(a\otimes\Delta_\ast(b))u_g$, where
$\Delta_ast$ is defined by $\Delta_\ast(b)(y)=b(\Delta(y))$ for all
$b\in\Lp^\infty(Y)$ and almost all $i\in Y$. It is easy to see that
two such endomorphisms $\varphi_\Delta,\varphi_\eta$ are unitarily
equivalent if an only if $\Delta=\eta$. The map
$\Delta\mapsto\varphi_\Delta$ embeds $\Factor(\Lambda\actson Y)^{\op}$
into $\End(M)$.

We give a set of conditions on the group $\Gamma_1$ and the action
$\Gamma\actson I$ that ensures that all endomorphisms
$\varphi:M\rightarrow M$ are of the form $\varphi_\Delta$ for some
$\Delta\in\Factor(\Lambda\actson Y)$. In fact, under the same
conditions, we also find the all right-finite $M$-$M$
bimodules are direct sums of bimodules of the form
$H_\Delta=\bimod{\varphi_\Delta(M)}{\Lp^2(M)}{M}$. In particular, we
see that $\Fundg_s(M)=\IN$.

In the next section, we give an explicit example of a group $\Gamma_1$
and an action $\Gamma\actson I$ that satisfy the conditions of theorem
\ref{thm:main}.

\begin{theorem}
  \label{thm:main}
  Let $\Lambda\actson (Y,\nu)$ be an ergodic probability measure
  preserving action. Assume that $\Lambda$ is anti-(T) as in
  definition \ref{def:anti-T} and that all 1-cocycles
  $\omega:\Lambda\times Y\rightarrow K$ with values in compact groups
  are trivial.

  Let $\Gamma_1$ ba a countable group that satisfies the conditions
  G$_1,\ldots, $G$_5$ below.
  \begin{itemize}
  \item[G$_1$] $\Gamma_1$ has the small normalizers property introduced above.
  \item[G$_2$] $\Gamma_1$ does not have non-trivial finite dimensional unitary representations.
  \item[G$_3$] all centralizers $\Centr_{\Gamma_1}\{g\}\subset\Gamma_1$ of finite-order elements
    $e\not=g\in\Gamma_1$ have the Haagerup property.
  \item[G$_4$] $\Gamma_1$ contains a property (T) subgroup $G$.
  \item[G$_5$] $G$ contains a non-amenable subgroup $G_0$ such that the
    commensurator of $G_0$ and $G$ generate all of $\Gamma_1$.
  \end{itemize}

  Let $\Sigma\subset \Gamma_1$ be an amenable subgroup that is not
  virtually abelian. Consider the amalgamated free product
  $\Gamma=\Gamma_1\free_\Sigma(\Sigma\times\Lambda)$. Let
  $\Gamma\actson I$ be an action of $\Gamma$ on a countable set,
  satisfying the conditions A$_1,\ldots,$A$_5$ below.
  \begin{itemize}
  \item[A$_1$] $\Gamma$ acts transitively on $I$
  \item[A$_2$] All stabilizers $\Stab\{i\}$ have the Haagerup property and
    $\Stab\{i\}\cap\Gamma_1$ is abelian.
  \item[A$_3$] the stabilizers $\Stab\{i,j\}$ of two-point sets are trivial.
  \item[A$_4$] The only injective group morphisms
    $\theta:\Gamma_1\rightarrow\Gamma_1$ that map
    each $\Stab_{\Gamma_1}\{i\}$ into some $\Stab_{\Gamma_1}\{j\}$ up
    to finite index are inner.
  \item[A$_5$] there is $i_0\in I$ such that $\Stab\{i_0\}\cap
    s\Gamma_1s^{-1}$ is infinite for all $s\in\Lambda$.
  \end{itemize}

  Let $X_0$ be an atomic probability space with unequal weights. Write
  $X=X_0^I$ and let $\Gamma$ act on $X$ by generalized Bernoulli
  action. Consider the natural quitient morphism
  $\pi:\Gamma\rightarrow\Lambda$ and 
  define a probability measure preserving action $\Gamma\actson X\times Y$ by
  by the formula $g(x,y)=(gx,\pi(g)y)$.
  Denote $M=\Lp^\infty(X\times Y)\rtimes\Gamma$.

  Then for every normal $\ast$-homomorphism
  $\varphi:M\rightarrow M$ from $M$ into itself is of the form
  \[\varphi((a\otimes b)u_g)=(a\otimes \Delta_\ast(b))u_g\text{ for
    some factor map }\Delta\in\Factor(\Lambda\actson Y),\]
  up to unitary conjugacy.

  Moreover, if $H$ is a right-finite dimensional $M$-$M$ bimodule,
  them $H$ is (isomorphic to) a finite direct sum of $M$-$M$ bimodules of the form
  $\bimod{\varphi_\Delta(M)}{\Lp^2(M)}{M}$ for factor maps
  $\Delta\in\Factor(\Lambda\actson Y)$.
\end{theorem}
\begin{proof}
  Observe that the statement about endomorphisms follows immediately
  from the statement about right-finite bimodules. So let $H$ be a
  right-finite bimodule of $M$. Then we know that $H$ is of the form
  $H=\bimod{\varphi(M)}{p(\IC^n\otimes \Lp^2(M))}{M}$ for some
  $\ast$-homomorphism $\varphi:M\rightarrow p(M\otimes\MatM_n(\IC))p$ from $M$ into
  some finite amplification $pM^np$ of $M$.

  We prove in 5 steps that $\Tr(p)=k$ for some natural number $k$
  and $\varphi$ is of the form
  \[\varphi((a\otimes b)u_g)=\sum_{i=1}^k e_{i,i}\otimes(a\otimes
  (\Delta_i)_\ast(b))u_g,\]
  up to unitary conjugacy, for some factor maps
  $\Delta_1,\ldots,\Delta_k\in\Factor(\Lambda\actson Y)$.

  Throughout the proof we will use the following notations for various
  subalgebras of $M$.
  \begin{itemize}
  \item $A=\Lp^\infty(X)$ and $B=\Lp^\infty(Y)$.    
  \item $M_1=(A\otimes B)\rtimes\Gamma_1$ and $M_2=(A\otimes
    B)\rtimes(\Sigma\times\Lambda)$.
  \item $M_{(i)}=(A\otimes B)\rtimes\Stab\{i\}$ for $i\in I$.
  \item $N=B\rtimes\Gamma$.
  \end{itemize}
  We will also combine the notations above. That way, we write
  $M_{1,(i)}$ for $M_1\cap M_{(i)}$. Similarly, we write $N_1=N\cap
  M_1$ and so on.

  We use the notation $M^n=M\otimes\MatM_n(\IC)$ for the amplification
  of $M$ by integer numbers. Similarly, we write
  $B^n=B\otimes\MatM_n(\IC)$ and $N^n=N\otimes\MatM_n(\IC)$. We denote
  the action of $\Gamma$ on $A$ by $\sigma$.

  \begin{step}
    We can assume that $p\in N_1^n$ and
    \[\varphi(M_1)\subset pM_1^np\text{ and
    }\varphi(N)\subset pN^np.\]
  \end{step}
  Since $\Gamma_2$ has the Haagerup property while $G$ has property
  (T), we see that $\varphi(\Lg(G))\nembeds_M M_2$.
  By \cite[theorem 5.1]{IoanaPetersonPopa:AFP}, we find a partial isometry $v\in M\otimes
  \MatM_{n,m}(\IC)$ with left support $vv^\ast=p$ and with right
  support $q=v^\ast v\in M_1^m$, and such that
  $v^\ast \varphi(\Lg(G))v\subset qM_1^mq$.
  We conjugate $\varphi$ by $v$ and assume already that $p\in
  M_1^n$ and $\varphi(\Lg(G))\subset
  pM_1^np$.

  Similarly, we see that $\varphi(\Lg(G))$ does not embed into
  $M_{(i)}$ for any $i\in I$. Applying \cite[Corollary
  4.3]{IoanaPopaVaes:vNsuperrigidity} (which is a version of
  \cite[theorem 4.1]{Popa:StrongRigidity1}) to the rigid inclusion
  $\varphi(\Lg(G))\subset pM_1^np$, we find a partial isometry $w\in
  M_1\otimes\MatM_{n,k}(\IC)$ with left support $p=ww^\ast$ and with
  right support $r=w^\ast w\in N_1^k$, satisfying
  $w\varphi(\Lg(G))w^\ast\subset rN_1^kr$. We conjugate $\varphi$ by $w$
  and we assume that $p\in N_1^n$ and that $\varphi(\Lg(G))\subset
  pN_1^np$.

  Observe that $\varphi(\Lg(G_0))$ is contained in $pN_1^np$, but it
  does not embed into $P$ nor into $N_{(i),1}$ for any $i\in I$. So
  by \cite[theorem 1.2.1]{IoanaPetersonPopa:AFP} and \cite[lemma
  4.2.1]{Vaes:Bimodules} (which is based on \cite[section
  3]{Popa:StrongRigidity1}), it follows that its quasi-normalizer
  is still contained in $pN_1^np$. But this quasi-normalizer contains
  $\varphi(B\otimes \Lg(\Comm_{\Gamma_1}(G_0)))$, and together with
  $\varphi(\Lg(G))$, this algebra generates $\varphi(N_1)$. We
  conclude that $\varphi(N_1)\subset pN_1^np$.

  We know that $\varphi(A)$ is an abelian subalgebra and all
  $\varphi(u_g)\in pM_1^np$ with $g\in G$ normalize
  $\varphi(A)$. Moreover, $\varphi(\Lg(G))$ does not embed into $P$. It
  is shown in \cite[theorem 1.4.1]{IoanaPetersonPopa:AFP} that then $\varphi(A)$ itself is contained in
  $pM_1^np$. Hence all of $\varphi(M_1)$ is contained in $pM_1^np$.

  We know that $\varphi(\Lg(\Lambda))$ commutes with
  $\varphi(\Lg(\Sigma))\subset pN^np$. Once we show that
  $\varphi(\Lg(\Sigma))$ does not embed into any $N_{(i)}$, inside $N$, then we can
  apply \cite[lemma 4.2.1]{Vaes:Bimodules} and conclude that $\varphi(\Lg(\Lambda))$ is
  contained in $pN^np$. In that case, we find that $\varphi(N)\subset
  pN^np$. It remains to show that $\varphi(\Lg(\Sigma))$ does not
  embed into $N_{(i)}$ for any $i\in I$. Observe that
  $N_{(i),1}=B\otimes\Lg(\Stab_{\Gamma_1}\{i\})$ is abelian, while
  $\Sigma$ is not virtually abelian. It follows that
  $\varphi(\Lg(\Sigma))$ does not embed into $N_{(i),1}$.
  Hence (see \cite[remark 3.3]{Vaes:Bimodules}) we find a sequence $(v_m)_m$ of unitaries in
  $\varphi(\Lg(\Sigma))$ such that
  \[\norm{\E_{N_{(i),1}}(xv_my)}_2\rightarrow 0\text{ for all }x,y\in
  N_1\text{ and for all }i\in I.\]
  We want to show that
  \[\norm{\E_{N_{(i)}}(xv_my)}_2\rightarrow 0\text{ for all }x,y\in
  N\text{ and for all }i\in I.\]
  By Kaplanski's density theorem, we can assume that $x=u_g,y=u_h$
  with $g,h\in\Gamma$. Write the fourier expansion of $v_m$ as
  $v_m=\sum_{k\in\Gamma_1}v_{k,m}u_k$. Then we compute that
  \begin{align*}
    \norm{\E_{N_{(i)}}(xv_my)}_2^2 &=
    \norm{\sum_{k\in\Gamma_1}\E_{N_{(i)}}(u_gv_{k,m}u_{kh})}_2^2\\
    &= \sum_{k\in\Gamma_1\cap
      g^{-1}\Stab\{i\}h^{-1}}\norm{v_{k,m}}_2^2.
  \end{align*}
  If this last sum is non-empty, then there is a $k_0\in \Gamma_1\cap
  g^{-1}\Stab\{i\}h^{-1}$. Then it follows that $\Gamma_1\cap
  g^{-1}\Stab\{i\}h^{-1}=\Stab_{\Gamma_1}\{g^{-1}i\}k_0$. So we see
  that
  \begin{align*}
    \norm{\E_{N_{(i)}}(xv_my)}_2^2
    &= \sum_{k\in\Stab_{\Gamma_1}\{g^{-1}i\}k_0}\norm{v_{k,m}}_2^2\\
    &= \norm{\E_{N_{(g^{-1}i),1}}(v_mu_{k_0}^\ast)}\rightarrow 0.
  \end{align*}

  We have shown that $\varphi(\Lg(\Lambda))$ does not embed into
  $N_{(i)}$ for any $i$. It follows that $\varphi(N)\subset
  pN^np$. This finishes the proof of our first step.

  \begin{step}
    We write $A_0=\Lp^\infty(X_0)$ and for all subsets
    $J\subset I$, we denote by
    $A_0^J=\Lp^\infty(X_0^J)$ the subalgebra of $A$ that consists of 
    functions that depend only on the components indexed by $J$.

    Let $J\subset I$ be an infinite subset that is invariant under an
    infinite group $H\subset\Gamma_1$.
    Then we show that $\varphi(A_0^{J})$ does not embed into $B$
    inside $M_1$, for any $i\in I$.
  \end{step}
  Observe that $\Gamma_1$ acts trivially on $B$, so $B$ is contained
  in the center of $M_1$, in fact, it is the center of $M_1$. Remark
  that any given element
  $e\not=g\in\Gamma_1$ can fix at most one $j\in I$. So $H$ acts freely
  on $A_0^J$, because $J$ is infinite.

  Now suppose that $\varphi(A_0^J)$ embeds into $B$ inside $M_1$. Then
  we find a non-zero partial isometry $v\in M_1\otimes \MatM_{n,m}(\IC)$ and a
  $\ast$-homomorphism $\theta:A_0^J\rightarrow qB^mq$, for some
  projection $q\in B^m$, such that $\varphi(x)v=v\theta(x)$ for all
  $x\in A_0^J$. But $B^m$ is of finite type I and $\theta(A_0^J)$ is
  a (non-unital) abelian subalgebra. Up to a unitary in $B^m$, we can
  assume that $q\in B\otimes\MatD_m(\IC)$ and $\theta(A_0^J)\subset
  q(B\otimes\MatD_m(\IC))q$ (see for example \cite[lemma C.2]{Vaes:Bourbaki} for
  an argument). Taking a non-zero component of $v$, we can assume that
  $m=1$ and hence that $v\in M_1\otimes\MatM_{n,1}(\IC)$. Set
  $r=vv^\ast\in \varphi(A_0^J)^\prime\cap pM_1^np$ and observe that
  $\varphi(x)r=r(1\otimes\theta(x))$ for all $x\in A_0^J$.

  Set $w_g=\varphi(u_g)$ for all $g\in H$.
  We claim that $r$ is orthogonal to $w_grw_g^\ast$ for all
  $e\not=g\in H$. Fix $g\in H$. Since $H$ acts freely on $A_0^J$, we
  find an element $a\in A_0^J$ such that $a-\sigma_g(a)$ is
  invertible. Then we compute that
  \[\varphi(a)rw_grw_g^\ast=r(1\otimes\theta(a))w_grw_g^\ast
  =rw_g(1\otimes\theta(a))rw_g^\ast=rw_gr\varphi(a)w_g^\ast=rw_grw_g^\ast\varphi(\sigma_g(a)).\]
  But we know that $rw_grw_g^\ast$ commutes with $\varphi(A_0^J)$, so
  $rw_grw_g^\ast\varphi(a-\sigma_g(a))=0$. Since $a-\sigma_g(a)$ is
  invertible, it follows that $r$ is orthogonal to $w_grw_g^\ast$.

  Since $H$ is an infinite group, we have found an infinite sequence
  of pairwise orthogonal projections with the same trace in the finite
  von Neumann algebra $pM_1^np$. This contradiction shows that
  $\varphi(A_0^J)$ can not embed into $B$ inside $M_1$.

  \begin{step}
    From now on, we allow $n=\infty$. In that case, we denote
    $\MatM_\infty(\IC)$ for $\Bounded(\ell^2(\IN))$, and we write
    $\MatM_{k,\infty}(\IC)$ for $\Bounded(\IC^k,\ell^2(\IC))$ and
    finally, we write $\MatD_\infty(\IC)=\ell^\infty(\IN)$.

    With these notations we can assume that $p\in B^n$ is a projection
    with finite trace and that
    $\varphi(A\otimes B)$ is
    contained in $p(A\otimes B^n)p$. Moreover, we can assume that
    $\varphi(u_g)=u_g$ for all $g\in\Gamma_1$.
  \end{step}
  Denote $C=\varphi(A)^\prime\cap pM_1^np$, and observe that
  $M_1=(A\rtimes\Gamma_1)\otimes B$, because $\Gamma_1$ acts trivially
  on $Y$.
  We want to apply theorem \ref{thm:inter-cartan} to conclude that $C$ embeds into
  $A\otimes B$ inside $M_1$. We check its four conditions. The first
  condition is satisfied by step 2 above. Observe
  that $\varphi(\Lg(G))$ has property (T) and hence can not embed into
  any of the amenable algebras $M_{(i),1}$ for any $i\in I$. This
  shows that the second and fourth condition are also satisfied.

  If the
  third condition were not satisfied, then we had that
  $\varphi(A\rtimes G)\embeds_{M_1} N_1$. So we find a partial
  isometry $0\not=v\in M_1\otimes\MatM_{n,m}(\IC)$ and a $\ast$-homomorphism
  $\theta:A\rtimes G\rightarrow qN_1^mq$, for some $q$, such that
  $\varphi(x)v=v\theta(x)$ for all $x\in A\rtimes G$. We can assume that $q$ is
  the support projection of $\E_{N_1}(v^\ast v)$. Observe that
  $\theta(A)$ is an abelian subalgebra of $qN_1^mq$ with large
  normalizer. Since $\Gamma_1$ has small normalizers, we see that
  $\theta(A)$ embeds into $B$ inside $N_1$. So there exists $0\not=w\in
  N_1\otimes\MatM_{m,k}(\IC)$ and a $\ast$-homomorphism
  $\rho:A\rightarrow rB^kr$ such that $\theta(x)w=w\rho(x)$. Since
  $vw\not=0$, it follows that $\varphi(A)$ embeds into $B$ inside
  $M_1$. But that is impossible by step 2 above. Hence also the third
  condition of theorem \ref{thm:inter-cartan} is satisfied.

  We conclude that $C\fembeds_{M_1} A\otimes B$. We want to apply
  theorem \ref{thm:inter-crit} to the inclusion $C\subset pM_1^np$. Denote by
  $\gamma:\Gamma_1\rightarrow\Unitary(pN_1^np)$ the group morphism
  that is defined by $\gamma(g)=\varphi(u_g)$. Observe that $\Gamma_1$
  does not have any non-trivial finite-dimensional representation and
  that $G\subset\Gamma_1$ is a property (T) subgroup. We show that
  $\gamma(\Gamma_1)^{\prime\prime}$ does not embed into
  $\Lg(\Centr_{\Gamma_1}\{g\})\otimes B$ for any
  $e\not=g\in\Gamma_1$. If $g$ is a finite-order element, then the
  centralizer has the Haagerup property, by assumption. If on the
  other hand $g$ has infinite order, then we know that
  $\Lg(g^{\IZ})\otimes B$ is a diffuse abelian subalgebra of
  $\Lg(\Gamma_1)\otimes B$. If $\gamma(\Gamma_1)^{\prime\prime}$
  embeds into the centralizer of $g$, then it follows that
  $\Lg(g^{\IZ})\otimes B$ has large normalizer. But $\Gamma_1$ is a
  group with small normalizers, so we conclude that
  $\Lg(g^{\IZ})\otimes B$ embeds into $B$, or still, that $g$ has
  finite order. This contradicts our assumption.
  We conclude that
  $\gamma(\Gamma_1)^{\prime\prime}$ does not embed into
  $\Lg(\Centr_{\Gamma_1}\{g\})\otimes B$ for any $e\not=g\in\Gamma_1$.

  We show that the action by conjugation on $\Centre(C)$ is weakly
  mixing relative to $D=\Centre(C)\cap pN_1^np$. Let
  $H\subset\Lp^2(\Centre(C))$ be a finite dimensional,
  $\gamma(\Lambda_1)$-invariant right $D$-submodule. Observe that then
  $HpN_1^np$ is a $\gamma(\Lambda_1)^{\prime\prime}$-$pN_1^np$ subbimodule of
  $\Lp^2(pM_1^np)$, and it has finite dimension on the right.
  Remark that $\gamma(\Gamma_1)^{\prime\prime}$ does not
  embed into $N_{(i),1}$ for any $i\in I$, so \cite[lemma 4.2.1]{Vaes:Bimodules} shows that
  any $\gamma(\Lambda_1)^{\prime\prime}$-$pN_1^np$ subbimodule of
  $\Lp^2(pM_1^np)$ that has finite dimension on the right, must be
  contained in $\Lp^2(pN_1^np)$. It follows that $H$ is contained in
  $\Lp^2(D)$.

  We have just shown that $C\fembeds_{M_1}A\otimes B$ and it is
  obvious from the definition of $C$ that $\Centre(C)^\prime\cap
  pM_1^np=C$, so we can apply theorem \ref{thm:inter-crit}. This yields a partial
  isometry $v\in N_1\otimes\Bounded(\IC^n,\ell^2(\IN)\otimes\ell^2(\IN))$ with left support
  $p=vv^\ast$ and with right support $q=v^\ast v\in
  \widetilde B=B\otimes\ell^\infty(\IN)\otimes\Bounded(\ell^2(\IN))$, and such that
  \[v^\ast Cv=A\otimes \widetilde B.\]
  Of course we still have that
  $v^\ast\varphi(N)v\subset q(N\otimes\Bounded(\ell^2(\IN))\otimes\Bounded(\ell^2(\IN)))q$.

  Moreover, we find a group morphism $\delta$ from $\Gamma_1$ to the
  group
  \[\cG=\left\{\left.\sum_{g\in\Gamma_1}p_gu_g\,\right\vert\, p_g\in
  B\otimes\ell^\infty(\IN)\text{ are projections with }\sum_g p_g=q\right\}.\]
  This group morphism satisfies $v^\ast\varphi(u_g)v=\delta(g)$ for all
  $g\in\Gamma_1$.

  Denote by $Z\subset Y\times\IN$ the
  support of $q$. Then we can view that group morphism
  $\Delta:\Gamma_1\rightarrow\cG$ as a measurable field
  $(\delta_z)_{z\in Z}$ of group morphisms
  $\delta_z:\Gamma_1\rightarrow\Gamma_1$.
  We consider $\delta(g)$ as a map from $Z$ to
  $\Lg(\Gamma_1)$, for all $g\in\Gamma_1$. As such,
  this map is given by $\delta(g)(z)=u_{\delta_z(g)}$.

  We first show that almost all the $\delta_z$ with $z\in Z$ are
  injective. Observe that $A\rtimes\Gamma_1$ is a factor and that
  $\varphi(A\rtimes\Gamma_1)$ is contained in
  $(A\rtimes\Gamma_1)\otimes q\widetilde Bq$. If not all $\delta_z$
  were injective, then there is an element $g\in\Gamma_1$ and a
  non-null set $U\subset Z$ such that $\delta_z(g)=e$ for all $z\in
  U$. Denote by $r=\chara_U q$ the central projection in
  $(A\rtimes\Gamma_1)\otimes q\widetilde Bq$ that corresponds to
  $U$. But then the $\ast$-homomorphism that maps $x\in
  A\rtimes\Gamma_1$ to $r\varphi(x)$ is not injective. This
  contradicts the factoriality of $A\rtimes\Gamma_1$.

  Fix $i\in I$ and suppose that there is a non-null set $V\subset Z$
  such that $\delta_z(\Stab_{\Gamma_1}\{i\})\cap\Stab\{j\}$  does not
  have finite index in $\delta_z(\Stab_{\Gamma_1}\{i\})$, for any
  $j\in I$ and for all $z\in V$. In other words,
  $\delta_z(\Stab_{\Gamma_1}\{i\})$ acts with infinite orbits on $I$,
  for all $z\in V$. Denote by $r_2=\chara_Vq$ the central projection in
  $(A\rtimes\Gamma_1)\otimes q\widetilde Bq$ that corresponds to
  $V\subset Z$. Then we see that $\varphi(A_0^{i})r_2\subset r_2\widetilde
  Br_2$. But we also get that
  $\varphi(A_0^{gi})r_2=\varphi(u_gA_0^{i}u_g)r_2$ is contained in
  $r_2\widetilde Br_2$ for all $g\in\Gamma_1$. Hence we get that
  $\varphi(A_0^{\Gamma_1i})\subset r_2\widetilde Br_2$, but this
  contradicts step 2. So we can conclude that for almost all $z\in U$
  and for every $i\in I$, a finite index subgroup of
  $\delta_z(\Stab_{\Gamma_1}\{i\})$ is contained in
  $\Stab_{\Gamma_1}\{j\}$ for some $j\in I$.

  We know that all such injective group morphisms from $\Gamma_1$ to
  itself are inner. So each $\delta_z$ is inner. This is the same this
  as saying that $\delta$ itself is conjugate to $g\mapsto 1\otimes
  u_g\in\cG$, inside $\cG$. So we find an element $u\in\cG$ such that
  $u\varphi(u_g)u^\ast=u_g$ for all $g\in \Gamma_1$. We conjugate
  $\varphi$ by $u$ and assume that $\varphi(u_g)=u_g$. Remark that $u$
  normalizes $\widetilde B$, so we still have that
  \[\varphi(A\otimes B)\subset A\otimes q\widetilde Bq
  \subset A\otimes
  q(B\otimes\Bounded(\ell^2(\IN))\otimes\Bounded(\ell^2(\IN)))q\cong
  A\otimes qB^\infty q.\]

  \begin{step}
    We can assume that there is a cocycle $(b_s)_{s\in\Lambda}$ with
    values in $\Unitary(pB^np)$ such that
    $\varphi(u_g)=b_{\pi(g)}u_g$ for all $g\in\Gamma$. Moreover, we
    can assume that
    $\varphi(a)=a$ for all $a\in A$.
  \end{step}
  Denote by $I_0\subset I$ the set of all $i\in I$ such that
  $\Stab_{\Gamma_1}\{i\}$ is infinite. Observe that, for each $i\in
  I_0$, we have that $\varphi(A_0^{\{i\}})$ commutes with
  $\Lg(\Stab_{\Gamma_1})$. But $\Lg(\Stab_{\Gamma_1})$ does not embed into
  $B\rtimes\Stab\{i,j\}$ for any $j\not=i$, simply because
  $\Stab\{i,j\}$ is trivial. Now we can apply \cite[lemma 4.2.1]{Vaes:Bimodules} and we
  obtain that $\varphi(A_0^{\{i\}})$ is contained in
  $A_0^{\{i\}}\otimes pB^np$.

  We show that
  $\varphi\restrict{A_0^{\{i\}}}$ is in fact the identity morphism on
  $A_0^{\{i\}}$, for all $i\in I_0$. Take an element $g\in\Gamma_1$ such that
  $gi\not=i$. Observe that $gi$ is still in $I_0$. Remember that $A_0=\Lp^\infty(X_0,\mu_0)$ and that
  $(X_0,\mu_0)$ is a purely atomic probability space with unequal
  weights. Take a minimal projection $q$ in $A_0$. Take a maximal
  abelian subalgebra $B_0\subset pB^np$ such that $q_1=\varphi(q)$ and
  $q_2=\varphi(\sigma_g(q))=\sigma_g(q)$ are both in $A\otimes
  B_0$. We know that $B_0\cong\Lp^\infty(Z,\eta)$ for some probability
  space $(Z,\eta)$. We consider $q_1$ and $q_2$ as measurable maps
  from $Z$ to $A$. Observe that $q_1(z)$ and $q_2(z)$ are independent
  for almost all $z\in Z$,
  i.e. $\tau(q_1(z)q_2(z))=\tau(q_1(z))\tau(q_2(z))=\tau(q_1(z))^2$. Write
  $f(z)=\tau(q_1(z))$ for all $z\in \widetilde Z$. Then we know that
  \begin{equation*}
    \int_{Z}f(z)^2 dz = \tau(q_1q_2) = \tau(q\sigma_g(q)) =
    \tau(q)^2=\left(\int_{Z}f(z) dz\right)^2
  \end{equation*}
  The only positive functions satisfying this condition are the
  constant functions, so we see that $f(z)=\tau(q)$ almost
  everywhere.

  We can do the same thing for all minimal projections
  $q_x=\chara_{\{x\}}^{\{i\}}$ in $A_0^{\{i\}}$. We find that
  $\varphi(q_x)(z)$ is a projection in $A_0^{\{i\}}$ with trace
  $\tau(\varphi(q_x)(z))=\mu_0(\{x\})$, for almost every $z\in Z$ and
  for all $x\in X_0$. Moreover, we know that
  $1=\sum_x q_x$, so we see that $1=\sum_x\varphi(q_x)(z)$ for almost
  all $z\in Z$. Since $\mu_0$ has unequal weights, it follows that
  $\varphi(q_x)(z)=q_x$ a.e, for all $x\in X_0$.
  
  We prove that, for every $s\in\Lambda$, we get $\varphi(u_s)=b_su_s$
  for some unitary $b_s\in pB^np$. Write
  $b_s=\varphi(u_s)u_s^\ast\in p(B^n\rtimes\Gamma)p$.
  By condition A$_5$, there is an element $i_0\in I$ such that
  $\Lambda i_0\subset I_0$. In particular, we see that $\Sigma i_0$
  and $s\Sigma i_0$ are contained in $I_0$. So $b_s$ commutes with
  $A_0^{\Sigma i_0}$. We write the Fourier expansion of $b_s$ as
  $b_s=\sum_{g\in\Gamma}b_{s,g}u_g$, where $b_{s,g}\in pB^np$.
  Since $b_s$ commutes with $A_0^{\sigma i_0}$, it follows that
  $a\otimes b_{s,g}=\sigma_g(a)\otimes b_{s,g}$ for all $g\in\Gamma$ and
  all $a\in A_0^{\Sigma i_0}$. If $b_{s,g}$ is nonzero, then we
  conclude that $a=\sigma_g(a)$ for all $a\in A_0^{\Sigma i_0}$. It
  follows that $g\in\Stab(\Sigma i_0)$. Since $\Sigma\subset \Gamma_1$
  is not abelian, we know that $\Sigma$ can not be contained in
  $\Stab_{\Gamma_1}\{i_0\}$. So $\Sigma i_0$ contains at least two
  elements and hence $g\in\Stab(\Sigma i_0)=\{e\}$. We can conclude
  that $b_s\in pB^np$, or still $\varphi(u_s)=b_su_s$ with $b_s\in\Unitary(pB^np)$.

  Since $\Gamma_1$ and $\Lambda$ generate the group $\Gamma$, we see
  that $\varphi(u_g)=b_{\pi(g)}u_g$ for all $g\in\Gamma$. It remains
  to show that $\varphi(a)=a$ for all $a\in A$. We know already that
  $\varphi(a)=a$ for all $a\in A_0^{\{i_0\}}$. Since $\Gamma$ acts
  transitively on $I$, the same holds for all $a\in A$.

  \begin{step}
    Conclude that theorem \ref{thm:main} holds.
  \end{step}
  Consider $p\in B^n$ as a map from $Y$ to $\MatM_n(\IC)$, or
  $\Bounded(\ell^2(\IN))$ if $n=\infty$. Then we see
  that $\Tr(p(y))$ is $\Lambda$-invariant. So by ergodicity of the
  $\Lambda$-action, it is constant, and up to conjugation by a unitary
  in $B^n$ we can assume that $p$ itself is constant, i.e. $p\in
  1\otimes\MatM_n(\IC)$. Reducing $n$ we can assume that $p=1$. Since
  $p$ was a finite projection, it follows that $n$ is now finite.
  Now we can consider $(b_s)_s$ as a cocycle for the action
  $\Lambda\actson Y$ with values in $\MatU_n(\IC)$. We assumed that
  all such cocycles are trivial, so up to conjugation with a unitary
  in $B^n$, we can assume that $b_s=1$ for all $s\in\Lambda$.

  Remark that $\varphi(B)\subset (A\otimes B^n)\cap
  (B^n\rtimes\Gamma)=B^n$ by steps 1 and 3. Since $B$ is abelian, we can assume that $\varphi(B)\subset
  B\otimes\MatD_n(\IC)$. Any such $\ast$-homomorphism is given by a
  quotient map $\Delta:Y\times\{1,\ldots,n\}\rightarrow Y$. This
  quotient map is $\Lambda$-equivariant because $\varphi(u_s)=u_s$ for
  all $s\in\Lambda$.
  The image $\Delta(Y\times\{k\})$ is $\Lambda$-invariant and
  non-null, so it must be all of $Y$ up to measure $0$. In other
  words, the formula $\Delta_k(y)=\Delta(y,k)$ defines a factor map
  for the action of $\Lambda$ in $Y$. This works for all
  $k=1,\ldots,n$, so we see that our original right-finite bimodule
  $H$ is a direct sum $H=\oplus_{k=1}^nH_\Delta$, finishing the proof
  of the theorem.
\end{proof}

\section{An example of a group action satisfying the conditions of
  theorem \ref{thm:main}}
\label{sect:ex}
In order to apply theorem \ref{thm:main}, we have to find an action
$\Gamma\actson I$ that satisfies the long list of conditions given
there. Such examples are necessarily rather complicated. This section
is devoted to the description of one such example.

As prescribed by theorem \ref{thm:main}, the group $\Gamma$ is of the form
$\Gamma=\Gamma_1\free_\Sigma(\Sigma\times\Lambda)$. We build the group
$\Gamma_1$ from an arithmetic lattice in $\Sp(n,1)$. We refer to
\cite{Morris:ArithmeticGroups} for an introduction to arithmetic lattices in Lie
groups. For this section, we only need to know that suitable
arithmetic subgroups are indeed lattices in the corresponding Lie groups.

Consider the set $\Hurwitz$ of Hurwitz quaternions, i.e.
\[\Hurwitz=\left\{a+bi+cj+dk\,\left\vert\, \text{ either }a,b,c,d\in \IZ\text{ or
}a,b,c,d\in\IZ+\frac12\right.\right\},\]
so the components are allowed to be either integers or half-integers,
but mixtures are not allowed. This is a ring under the usual
addition and multiplication of quaternions, so
$i^2=j^2=k^2=-1$ and $ij=k=-ji$. We denote the element
$\frac12(1+i+j+k)$ by $h$. The skew field of quaternions is denoted by $\IH$.

The quaternions come with a natural involution defined by
$\overline{a+bi+ci+dk}=a-bi-ci-dk$. This involution reverses the order
of the multiplication,
i.e. $\overline{x\,y}=\overline{x}\,\overline{y}$.
Consider the sesquilinear form
$B:\Hurwitz^{3}\times\Hurwitz^{3}\rightarrow\Hurwitz$ on
$\Hurwitz^{3}$ that is defined by
$B(\xi,\eta)=\overline{\xi_0}\eta_0-\sum_{i=1}^2\overline{\xi_i}\eta_i$. Observe
that this form is of signature $(2,1)$. Consider the group
\[G=\PSp(B,\Hurwitz)=\{A\in\MatM_{3}(\Hurwitz)\mid
B(A\xi,A\eta)=B(\xi,\eta)\text{ for all
}\xi,\eta\in\Hurwitz^{3}\}/\{\pm 1\}.\]
This group is an arithmetic lattice in the Lie group $\Sp(2,1)/\{\pm1\}$. As
such it has property (T), and by \cite{PopaVaes:CartanHyperbolic}, $G$ has the small
normalizers property.

We remark that $\SL_2\IZ$ embeds into $G$, in the
following way. Observe that $\SL_2\IZ$ is exactly the set of
matrices in $\SL_2\IZ[i]$ that preserve the non-definite Hermitian
form
$B_0(\xi,\eta)=\overline{\xi_2}i\eta_1-\overline{\xi_1}i\eta_2$. Consider
the linear transformation $A:\IZ[i]^2\rightarrow\Hurwitz^2\subset\Hurwitz^3$
that is defined by the matrix
\[A=\left(\begin{matrix}h&\overline{h}i\\i\overline{h}&-ihi\end{matrix}\right).\]
Observe that $B(A\xi,A\eta)=B_0(\xi,\eta)$, and that the matrix $A$ is
invertible over $\Hurwitz$. So $A$ defines an embedding of $\SL_2\IZ$
into $G$ mapping a matrix $B\in\SL_2\IZ$ to the block matrix
\[\left(\begin{matrix}ABA^{-1}&0\\0&1\end{matrix}\right)\].

Consider the subgroup $G_0$ of all elements in $\SL_2\IZ$ that
are represented by matrices of the form
\[\left(\begin{matrix}a&b\\c&d\end{matrix}\right)\qquad
  \begin{aligned}
    &\text{where }a,b,c,d\text{ are integers}\\
    &\text{ and }b=0\mod 6, a=d=1\mod 6\text{ and } c=0\mod 7.\end{aligned}\]
Define the group $\Gamma_1=G\free_{G_0}\SL_2\IQ$, and
consider its subgroup $\Sigma=\ST_2\IQ\subset\SL_2\IQ$ of upper
triangular matrices in $\SL_2\IQ$. Now it
follows from theorem \ref{thm:smallnorm-afp} that $\Gamma_1$ has the
small normalizers property. The group $G\subset\Gamma_1$ has property
(T). We also see that $G_0$ is almost
normal in $\SL_2\IQ$. But $G$ and $\SL_2\IQ$ together generate all
of $\Gamma_1$. It was shown in \cite{vNW:MinimallyAlmostPeriodic} that $\SL_2\IQ$ does not
have any non-trivial finite-dimensional representations. Because the
group $G_0$ generates all of $G$ as a normal subgroup, the same
is true for $\Gamma_1$. Moreover, the group $\Sigma$ is amenable, but
not virtually abelian.

We show that $\Centr_{\Gamma_1}\{g\}$ has the Haagerup property for
every finite-order element $e\not=g\in\Gamma_1$. Observe that
a finite-order element $e\not=g$ in $\Gamma_1$ is conjugate to an
element in one of the components $G$ or $\SL_2\IQ$. We can assume
without loss of generality that $g\in G$ or $g\in \SL_2\IQ$. Moreover, $G_0$ is
torsion-free, so $g$ is not conjugate to an element in $G_0$. Hence
the centralizer $\Centr_{\Gamma_1}\{g\}$ is still contained in the same component $G$
respectively $\SL_2\IQ$.
If $g$ were in $\SL_2\IQ$, then this already implies that the
centralizer of $g$ has the Haagerup property. On the other hand, if
$g$ were in $G$, then the centralizer in $\Gamma_1$ is just the
centralizer in $G$. This last centralizer is a discrete subgroup of
the centralizer $C$ of $g$ in the Lie group $\Sp(2,1)/\{\pm1\}$. The centralizer $C$
in $\Sp(2,1)/\{\pm1\}$ is a Lie group of strictly smaller dimension than $\Sp(2,1)/\{\pm1\}$, hence its
Lie algebra does not contain a copy of $\liesp(2,1)$ nor of
$\liesl_2\IR\ltimes\IR^2$. Now it follows from \cite{CCJJV:Haagerup} or
\cite{deCornulier:HaagerupSL2} that $C$ has the Haagerup property. The same is
true for its discrete subgroup $\Centr_{\Gamma_1}\{g\}$.

Till now, we have checked all the conditions on the group $\Gamma_1$
that do not depend on the action $\Gamma\actson I$. We define the action
$\Gamma\actson I$ as follows. Choose a one-to-one map
$\Lambda\ni\lambda\mapsto n_\lambda\in\IN$. For every
$\Lambda\in\Lambda$, we set $x_\lambda=n_\lambda+i$. Consider the matrix
\[B_\lambda=\left(\begin{matrix}\overline
  x_\lambda&0&n_\lambda\\0&1&0\\n_\lambda&0&x_\lambda\end{matrix}\right)\in
  \Hurwitz^{3\times 3}.\]
This matrix defines an element in $G$. Consider likewise the element
\[C_\lambda=\left(\begin{matrix}\frac1{n_\lambda}&1\\-1&0\end{matrix}\right)\in\SL_2{\IQ}.\]
Finally, we define the element $h_\lambda=B_\lambda C_\lambda\in
\Gamma_1$ and we consider the subgroup $H\subset \Gamma$ generated by
the $\lambda h_\lambda\lambda^{-1}$. Finally, we define $I$ to be the
set of left cosets of $H$, i.e $I=\Gamma/H$, with the natural action
of $\Gamma$ by left translation.

Observe that all stabilizers of this
action are conjugate to $H$, which in turn is isomorphic to
$\FG_\infty$. In particular, the stabilizers all have the Haagerup
property. Moreover, the intersection of any conjugate of $H$ with
$\Gamma_1$ is either trivial or a copy of $\IZ$. So the stabilizers
$\Stab_{\Gamma_1}\{i\}$ are abelian. Remark also that
$H=\Stab_{\Gamma}\{e\}$ has an infinite intersection with all $\lambda
\Gamma_1\lambda^{-1}$, so $\Stab_{\Gamma_1}\{\lambda H\}$ is infinite
for all $\lambda\in\Lambda$. 

It remains to check that all stabilizers of two-point sets are
trivial, and that injective group morphisms
$\delta:\Gamma_1\rightarrow\Gamma_1$ that map stabilizers into
stabilizers up to finite index are inner. We prove a lemma that
implies both facts.
\begin{lemma}
  Let $e\not=a,b\in H$ be elements of $H$, let $g$ be an element of $\Gamma$
  and consider an injective group morphism
  $\delta:\Gamma_1\rightarrow\Gamma_1$. Observe that $\delta$ defines
  a group morphism (we abuse the notation and keep using the letter
  $\delta$) $\delta:H\rightarrow\Gamma$ by the formula $\delta(\lambda
  h_\lambda \lambda^{-1})=\lambda\delta(h_\lambda)\lambda^{-1}$. If we have
  $a=g\delta(b)g^{-1}$, then it follows that $\delta$ is inner,
  say $\delta=\Ad_h$ and moreover $gh\in H$.
\end{lemma}
\begin{proof}
  We begin by studying the injective group morphisms
  $\delta:\Gamma_1\rightarrow\Gamma_1$.
  \begin{step}
    Every injective group morphism
    $\delta:\Gamma_1\rightarrow\Gamma_1$ is bijective and moreover it
    is given by
    $\delta=\Ad_h\circ(\Ad_{E}\free_{G_0} \Ad_F)$ where $h\in\Gamma_1$ and either
    \begin{alignat*}{3}
        F&=\left(\begin{matrix}1&0\\0&1\end{matrix}\right)
        &\text{ and }
        E&=\left(\begin{matrix}V&0\\0&u\end{matrix}\right)
        &\text{ with }
        V&=A\left(\begin{matrix}v&0\\0&v\end{matrix}\right)A^{-1}\\
%
        \text{or }F&=\left(\begin{matrix}1&0\\0&-1\end{matrix}\right)
        &\text{and }
        E&=\left(\begin{matrix}V&0\\0&u\end{matrix}\right)
        &\text{ with }
        V&=A\left(\begin{matrix}vk&0\\0&-vk\end{matrix}\right)A^{-1}\\
    \end{alignat*}
    where $u\in \cU=\left\{\left.\pm \alpha, \frac{\pm \beta\pm \gamma}{\sqrt2}, \frac{\pm1\pm i\pm j\pm k}{2}\,\right\vert\, \alpha,\beta,\gamma=1,i,j,k\right\}$ and $v\in\cU\cap \IC$.
  \end{step}
  Let $\delta:\Gamma_1\rightarrow\Gamma_1$ be an injective group
  morphism. Since $\delta(G)$ has property (T), it must be contained
  in a conjugate of one of the two components $G$, $\SL_2\IQ$ of
  $\Gamma_1$. It can not be contained in a conjugate of $\SL_2\IQ$
  because that group has the Haagerup property. So we find
  $y_0\in\Gamma_1$ such that $y_0\delta(G)y_0\subset
  G$.
  
  Consider Zariski the closure of $G_1$ of $y_0\delta(G)y_0^{-1}$ inside the
  Lie group $\Sp(2,1)/\{\pm1\}$. Then we know that $G_1$ is a Lie
  group. If $G_1$ were not equal to $\Sp(2,1)/\{\pm1\}$, then we
  know that its Lie algebra is also strictly smaller that
  $\liesp(2,1)$. So it does not contain a Lie subalgebra of the form
  $\liesp(2,1)$ nor of the form $\liesl_2\IR\ltimes\IR^2$. By \cite{CCJJV:Haagerup,
    deCornulier:HaagerupSL2}, it follows that $G_1$ has the Haagerup property,
  which is absurd because it contains the discrete property (T)
  group $y_0\delta(G)y_0^{-1}$. So we conclude that $G_1=\Sp(2,1)/\{\pm 1\}$.
  Now the Margulis superrigidity theorem (or better said, its
  version for $\Sp(n,m)$, see \cite{Corlette:MargulisSuperrigidity}) shows that
  $\Ad_{y_0}\circ\delta\restrict{G}$ extends to an isomorphism of
  $\Sp(2,1)/\{\pm1\}$. All such isomorphisms are inner, so we find an
  element $y$ in the product $(\Sp(2,1)/\{\pm1\})\,\Gamma_1$ such that
  $\delta(x)=y^{-1}xy$ for all $x\in G$.

  On the other hand, we know that the element
  \[1+e_{1,2}=\left(\begin{matrix}1&1\\0&1\end{matrix}\right)\in\SL_2\IQ\]
  has roots of all order. The only elements in $\Gamma_1$ that have
  roots of all orders are conjugate to $1+re_{1,2}$ for some
  $r\in\IQ$. Hence we find an element $z\in \GL_2\IQ\Gamma_1$ such
  that $\delta(1+e_{1,2})=z^{-1}(1+e_{1,2})z$. Write
  \[x_v=\left(\begin{matrix}v&0\\0&v^{-1}\end{matrix}\right)\in\SL_2\IQ\text{ with
  }v\in\IR\setminus 0.\]
  Then we know that $x_v(1+e_{1,2})x_v^{-1}=1+v^2e_{1,2}$. It follows
  that $\delta(x_v)=z^{-1}x_vz$ and then it also follows that
  $\delta(1+e_{2,1})=z^{-1}(1+e_{2,1})z$. Since the $1+re_{1,2}$ and
  $1+re_{2,1}$ generate $\SL_2\IQ$, it follows that
  $\delta(x)=z^{-1}xz$ for all $x\in \SL_2\IQ$.

  For every element $x\in G_0$, we must have $\delta(x)=y^{-1}xy$ and
  $\delta(x)=z^{-1}xz$. It follows that $zy^{-1}$ commutes with
  $G_0$. The only elements in $\GL_2\IQ\Gamma_1(\Sp(1,2)/\{\pm1\})$
  that commute with $G_0$ are of the form $EF^{-1}$ where $E$ and $F$ are
  of the form described earlier. It follows that $\delta$ is indeed of
  the required form.

  Form now on, we assume that $\delta$ is of the form
  $\Ad_E\free\Ad_F$ where $E$ and $F$ are as in step 1.
  \begin{step}
    Replacing $g$ by an element in $Hg\delta(H)$ and adapting $a,b$
    accordingly, we can assume that $g$ is of the form $\lambda
    g_1\mu^{-1}$, while $a=\lambda h_{\lambda}^{s}\lambda^{-1}$ and
    $b=\mu \delta(h_{\mu})^t\mu^{-1}$ for some $g_1\in\Gamma_1$, some
    $\lambda,\mu\in\Lambda$ and $s,t\in\IZ$.
  \end{step}
  Among the elements in $Hg\delta(H)$, we can assume that $g$ has
  minimal length. Write $g=\lambda_0g_1\lambda_1\ldots g_n\lambda_n$
  where $\lambda_0,\lambda_1\in\Lambda$,
  $\lambda_1,\ldots,\lambda_{n-1}\in\Lambda\setminusb\{e\}$ and
  $g_1,\ldots,g_n\in \Gamma_1\setminus \Sigma$. Write $b$ as a reduced
  word $b=\mu_1 b_1\mu_1^{-1}\mu_2
  b_2\mu_2^{-1}\ldots\mu_mb_m\mu_m^{-1}$ where $b_i\in
  h_{\mu_i}^{\IZ}$ for all $i$. Then we study 5 cases.

  \textit{case 1:} $n\geq 2$, $m\geq 2$.\\ 
  In that case, we know that
  $g_1\not\in h_{\lambda_0}^{\IZ}\Sigma$ and that
  $g_n\not\in\Sigma\delta(h_{\lambda_n^{-1}})^{\IZ}$. The word
  \[g\delta(b)g^{-1}=\lambda_0g_1\lambda_1\ldots g_n\lambda_n \mu_1
  \delta(b_1)\mu_1^{-1}\ldots\mu_m\delta(b_m)\mu_m^{-1}\lambda_n^{-1}g_n^{-1}\ldots
  g_1^{-1}\lambda_0^{-1}\]
  is reduced, except that maybe $\lambda_n=\mu_1^{-1}$ or
  $\lambda_n=\mu_m^{-1}$. But even then, we know that $g_n\delta(b_1)$
  and $\delta(b_m)g_n^{-1}$ can not be contained in $\Sigma$. In any
  case, we get a reduced word for $g\delta(b)g^{-1}$ that starts with
  something of the form $\lambda_0g_1$ where $g_1\not\in
  h_{\lambda_0}^{\IZ}\Sigma$. Such an element can never be contained
  in $H$.

  \textit{case 2:} $n=1$, $m\geq 2$ and $g_1\not\in
  h_{\lambda_0}^{\IZ}\Sigma\delta(h_{\lambda_1})^{\IZ}$\\
  In this case, the argument is the same as for the first case: the
  obvious word for $g\delta(b)g^{-1}$ is almost reduced, and it can
  never be contained in $H$.

  \textit{case 3:} $n=1$ and $g_1\in
  h_{\lambda_0}^{\IZ}\Sigma\delta(h_{\lambda_1})^{\IZ}$.\\
  We can assume that $g_1\in\Sigma$, so $g$ is of the form $\lambda
  g_1$. For $g\delta(b)g^{-1}$ to be in $H$, we need at least that
  $\delta(b_m)g_1^{-1}\in \Sigma h_{\lambda\mu_m}^{\IZ}$. A direct
  computation shows that this is only possible if $g_1=e$. 
  So we have that $\delta(b_m)=k_1h_{\lambda\mu_1}^s$ for some
  $k_1\in\Sigma$ and $s\in\IZ$. But then it follows that $\delta(b_{m-1})k_1^{-1}$ is
  contained in $\Sigma h_{\lambda\mu_{m-1}}^{\IZ}$. As before, it follows
  that $k_1=e$. By induction, we see that $\delta(b_i)\in
  h_{\lambda\mu_i}^{\IZ}$ for all $i$. Replacing $b$ by any $\mu_i
  b_i\mu_i^{-1}$, we can assume that $b$ is of the required form. The
  conjugating element $g$ is already of the required form, and $a$ is
  then automatically of the right form.

  \textit{case 4:} $n\geq 2$, $m=1$\\
  In this case, we can assume that $g_1\not\in
  h_{\lambda_0}^{\IZ}\Sigma$.
  We get the following word for $g\delta(b)g^{-1}$:
  \[g\delta(b)g^{-1}\lambda_0 g_1\ldots \lambda_{n-1}g_n \lambda_n\mu_1
  \delta(b_1)\mu^{-1}\lambda_n^{-1}g_n^{-1}\lambda_{n-1}^{-1}\ldots
  g_1^{-1}\lambda_0^{-1}.\]
  If $\lambda_n$ is not $\mu_1^{-1}$, this word is already reduced. On
  the other hand,
  form the definition of $H$ and $\Sigma$, it is clear that no element
  of the form $\delta(h_{\lambda})^s$ can be conjugated into $\Sigma$,
  inside $\Gamma_1$. So if $\lambda_n\mu_1$ were $e$, then we still
  had that
  \[g\delta(b)g^{-1}\lambda_0g_1\ldots \lambda_{n-1} (g_n\delta(b_1)g_n^{-1}) \lambda_{n-1}^{-1} \ldots g_1^{-1}\lambda_0^{-1}\]
  is a reduced word. In both cases, it is clear that
  $g\delta(b)g^{-1}$ can not be contained in $H$.

  \textit{case 5:} $n=1$, $m=1$ and $g_1\not\in
  h_{\lambda_0}^{\IZ}\Sigma\delta(h_{\lambda_1})^{\IZ}$\\
  Now we see that
  \[a=g\delta(b)g^{-1}=\lambda_0
  g_1\lambda_1\mu_1\delta(b_1)\mu_1^{-1}\lambda_1^{-1}g_1^{-1}\lambda_0^{-1}.\]
  If $\lambda_1\mu_1\not=e$, then this is a reduced expression
  for an element that is not in $H$. Otherwise, $a$, $b$ and $g$ are
  of the required form. This finishes the proof of step 2.

  \begin{step}
    It follows that $\delta=\id$, $s=t$, $\mu=\lambda$ and
    $g=h_{\lambda}^r$ for some $r\in\IZ$.
  \end{step}
  We know that $g=\lambda g_1\mu^{-1}$, $a=\lambda h_\lambda^s\lambda^{-1}$
  and $b=\mu h_\mu^t\mu^{-1}$, and we also have that
  $a=g\delta(b)g^{-1}$. It follows that $g_1\delta(h_\mu^t)g_1^{-1}=h_{\lambda}^s$

  Observe that $h_{\lambda}^s=B_\lambda C_\lambda\ldots B_\lambda
  C_\lambda$ is a reduced word for $h_\lambda^s$. Moreover,
  $h_\lambda^s\in\Gamma_1$ has minimal length among its conjugates, in
  the amalgamated free product decomposition of $\Gamma_1$. The same
  is true for $\delta(h_\mu)$. It follows that $s=\pm t$. We can assume
  that $g_1$ has minimal length among the elements
  $h_\lambda^{r_1}g_1\delta(h_\mu)^{r_2}$ with $r_1,r_2\in\IZ$. With
  this assumption, it also follows that $g_1\in \{e,C_\lambda^{-1},B_\lambda\}G_0\{e,\delta(C_\mu),\delta(B_\mu)^{-1}\}$.

  In any case, it follows that there are $g_2,g_3\in G_0$ such that
  $g_2\delta(B_\mu)^{\sigma_2}=B_\lambda^{\sigma_3}g_3$ where
  $\sigma_2,\sigma_3=\pm 1$. Moreover, we get that $g_1=bg_2c$ where
  $b\in\{e,C_\lambda^{-1},B_\lambda\}$ and
  $c\in\{e,\delta(C_\lambda),\delta(B_\lambda)^{-1}\}$. We also see
  that $\sigma_2\sigma_3=-1$ if and only if either $b$ or $c$ is $e$,
  but not both.

  Using the fact that
  $\delta=\Ad_E\free\Ad_F$ with $E,F$ as in step 1, we see that
  $g_2EB_\mu^{\sigma_2}=B_\lambda^{\sigma_3} g_3E$. We write $g_2$ and $g_3$ as block
  matrices
  \[g_2=\left(\begin{matrix}AXA^{-1}&0\\0&1\end{matrix}\right)
    \text{ and }
    g_3=\left(\begin{matrix}AYA^{-1}&0\\0&1\end{matrix}\right),\]
  where $X$ and $Y$ are $2\times 2$ matrices with integer coefficients
  that are upper triangular with eigenvalues $1$, modulo $6$, and that
  are lowe triangular with arbitrary eigenvalues modulo 7. Denote
  $\tau^{\sigma_2}(x_\mu)$ for $x_\mu$ if $\sigma_2=1$ and
  $\overline x_\mu$ if $\sigma_2=-1$. We use also the similar
  notations with $\sigma_3$ and $x_\lambda$. Then we
  find that
  \begin{align}
    AXVA^{-1}{\scriptstyle\left(\begin{matrix}\overline {\tau^{\sigma_2}(x_\mu)}&0\\0&0\end{matrix}\right)}
    &=\sigma_1{\scriptstyle\left(\begin{matrix}\overline {\tau^{\sigma_3}(x_\lambda)}&0\\0&0\end{matrix}\right)}AYVA^{-1}
    \label{eqn:ex:step3:ul}\\
    AXVA^{-1}{\scriptstyle\left(\begin{matrix}n_\mu\\0\end{matrix}\right)}
    &=\sigma_1\sigma_2\sigma_3{\scriptstyle\left(\begin{matrix}n_\lambda\\0\end{matrix}\right)}u
    \label{eqn:ex:step3:ur}\\
    u{\scriptstyle\left(\begin{matrix}n_\mu\\0\end{matrix}\right)}
    &=\sigma_1\sigma_2\sigma_3{\scriptstyle\left(\begin{matrix}n_\lambda\\0\end{matrix}\right)}AYVA^{-1}
    \label{eqn:ex:step3:ll}\\
    \text{and }u\tau^{\sigma_2}(x_\mu)&=\sigma_1 \tau^{\sigma_3}(x_\lambda) u
    \label{eqn:ex:step3:lr}
  \end{align}
  where $V$ and $u$ are as in the definition of $E$ in step 1, and
  $\sigma_1=\pm 1$.
  The equation (\ref{eqn:ex:step3:lr}) implies that
  $\abs{x_\mu}=\abs{x_\lambda}$, so it follows that
  $\mu=\lambda$. Moreover, comparing the real parts in
  (\ref{eqn:ex:step3:lr}), we see that the sign $\sigma_1=1$.
 
  From equation (\ref{eqn:ex:step3:ur}), we conclude that
  \[\xi=A^{-1}{\scriptstyle\left(\begin{matrix}1\\0\end{matrix}\right)}={\scriptstyle\left(\begin{matrix}h\\-i\overline{h}\end{matrix}\right)}\]
  is an eigenvector of $XV$ with eigenvalue $x=\sigma_2\sigma_3u$.

  We have to take a little care with what we
  mean by eigenvalues and eigenvectors over non-abelian rings, but the
  definition is exactely as in (\ref{eqn:ex:step3:ur}): the matrix is
  on the left of the vector, while the scalar is on the right.
  In a skew field, like $\IH$, we say that two elements $x,y$ are
  equivalent if they are conjugates of each other. Observe that in
  $\IH$ two elements are equivalent if and only if they have the same
  modulus and the same real part. For the vector spaces $\IH^n$, we
  always consider scalar multiplication on the right and matrix
  multiplication on the left. If $Z$ is a square matrix over $\IH$ and
  $\eta$ is an eigenvector of $Z$ with eigenvalue $y\in\IH$. Then we
  see that a scalar multiple $\xi z$ is still an eigenvector of $Z$
  but with eigenvalue $z^{-1}yz$.

  Now we see that the vector $(1,ih)$ is an eigenvector of $XV$ with
  eigenvalue $\tilde u=\sigma_2\sigma_3hu\overline h$. Observe that
  $\tilde u$ is still in $\cU$. Write $X$ and $V$ as matrices
  \[X=\left(\begin{matrix}a&b\\c&d\end{matrix}\right)\text{ and }
    V=\left(\begin{matrix}v&0\\0&v\end{matrix}\right)\text{ or }
    V=\left(\begin{matrix}vk&0\\0&-vk\end{matrix}\right),\]
  where $v\in\cU\cap\IC$ and $a,b,c,d$ are integers.
  Because $(1,ih)$ is an eigenvector, we see that
  $cv+dvih=ih(av+bvih)$ or $cvk-dvkih=ih(avk-bvkih)$.
  Observe that $v, vih, ihv$ and $ihvih$ are linearly independent over
  $\IR$ unless $v=\pm1$. We know that $X\not=0$, so we conclude that
  $v=\pm 1$. Similarly, $vk, vkih, ihvk$ and $ihvkih$ are
  linearly independent over $\IR$ unless $v=\pm\frac{1-i}{\sqrt{2}}$.

  We check both cases separately. Suppose we are in the first case,
  with $v=\pm1$.One easily computes that $c=-b$ and $a=d+b$. Since $X\in G_0$, we
  see that $b$ is a multiple of $6$ and $d$ is congruent to $1$ modulo
  $6$. Remark that $\tilde u=av+bvih$, while $b$ is a multiple of $6$ and
  $a$ is congruent to $1$ modulo $6$. The only element in $\cU$ that
  can be written in this form is $v$ itself, so we conclude that $\tilde
  u=v$, i.e $u=\sigma_2\sigma_3v=\pm 1$. But if $\sigma_2\sigma_3=-1$,
  then equation (\ref{eqn:ex:step3:lr}) tells us that
  $x_\lambda=\overline{x_\lambda}$, which is simply not the case. So
  $u=v$ and hence $\delta=\id$. Moreover, we see that
  $X=1$ and so $g_1=1$. So in the first case, we are done.

  Assume now that we are in the second case and
  $v=\pm\frac{1-i}{\sqrt{2}}$. A direct computation shows that $a=d$
  and $c=-b-a$. Moreover, $a,-b,d$ are congruent to $1$ modulo $6$
  while $c$ is a multiple of $6$, and $b$ is a multiple of $7$. Since $avk-bvkih=\tilde u\in\cU$, we
  see that this is not possible. This contradiction shows that we must
  have been in the first case, and the lemma is proven.
\end{proof}

\section{Examples of endomorphism semigroups}
\label{sect:semigroups}
In this section, we use the results of sections \ref{sect:main} and
\ref{sect:ex} in order to show that many semigroups appear as
$\End(M)$ for some type II$_1$ factor $M$. The starting point is the following.
\begin{theorem}
  \label{thm:ex-finite}
  Let $\Lambda\actson (Y,\nu)$ be any probability measure preserving
  action of a not necessarily discrete group. Then there is a type II$_1$ factor $M$ such that
  \begin{align*}
    \End(M)&\cong\Factor(\Lambda\actson (Y,\nu))^{\op}\\
    \text{and }\RFBimod(M)&\cong\{\text{formal finite direct sums of elements of
    }\Factor(\Lambda\actson (Y,\nu))^{\op}\}
  \end{align*}
\end{theorem}
\begin{proof}
  First of all, we can assume that $\Lambda\subset
  \Autmp(Y,\nu)$. Take any countable dense subgroup $\Lambda_0$ of
  $\Lambda$. We observe that $\Factor(\Lambda\actson
  Y)=\Factor(\Lambda_0\actson Y)$. So we can assume that $\Lambda$ is
  a countable group.

  Then this theorem is a direct consequence of our main result
  \ref{thm:main}, the flexible class of examples in section
  \ref{sect:ex}, and of lemma \ref{lem:no-fd} below.
\end{proof}

Before we prove lemma \ref{lem:no-fd}, we show that many left
cancellative semigroups appear this way. First, observe that any
\emph{compact} left cancellative semigroup with unit is in fact a
group. But then we get that $G=\Factor(G\actson (G,h))$ where $h$
denotes the Haar measure on $G$. So the compact case is easy to
handle, but nor very interesting.

A more interesting class of left cancellative semigroups is the class
of discrete left cancellative semigroups. This class contains proper
semigroups, like $\IN$, and even semigroups that can not be embedded
into groups (for example they are not right cancellative). We show
that any discrete left cancellative semigroup appears as the (opposite
of the) semigroup
of factors of some probability measure preserving action.

Let $G$ be a left cancellative semigroup. Let $(Y_0,\nu_0)$ be a
standard nonatomic probability space and consider
$(Y,\nu)=(Y_0,\nu_0)^G$. Every element $g\in G$ defines a measure
preserving quotient map $f_g:Y\rightarrow Y$ by the formula
$f_g(x)_h=x_{gh}$. Then we see that $f_gf_h=f_{hg}$. So we see that
$G^{\op}\subset\Factor(Y,\nu)$. If we take a non-atomic base space
$(Y_0,\nu_0)$, and we choose an appropriate subgroup $\Lambda\subset
\Autmp(Y,\nu)$, then we get that $G^{op}=\Factor(Y,\nu)$:
\begin{lemma}
  \label{lem:discrete-factor}
  Let $G$ be a left cancellative semigroup with unit $e$, and let $(Y_0,\nu_0)$ be a
  non-atomic probability space. Set $(Y,\nu)=(Y_0,\nu_0)^G$ as before.
  Then there is a subgroup $\Lambda\subset\Autmp(Y,\nu)$ such that
  \[\Factor(Y,\nu)=G^{op}.\]
\end{lemma}
\begin{proof}
  Denote by $\cG\subset \Autmp(Y,\nu)$ the closed subgroup of
  $\Autmp(Y,\nu)$ that is generated by the following transformations:
  \begin{align*}
    \psi_\Delta&:Y\rightarrow Y&\psi_\Delta(x)_h&=\Delta(x_h)&&\text{for all }\Delta\in\Autmp(Y_0,\nu_0)\\
    \varphi_{U,g,\Delta}&:Y\rightarrow Y&\psi_\Delta(x)_h&=\begin{cases}\Delta(x_h)&\text{if }x_{hg}\in U\\x_h&\text{if }x_{hg}\not\in U\end{cases}&&\begin{aligned}&\text{for all }\Delta\in\Autmp(Y_0,\nu_0)\text{ with }\Delta(U)=U\\&\text{ and for all }g\in G, U\subset Y_0\end{aligned}\\
  \end{align*}
  These automorphisms commute with all the $f_g$ with $g\in G$. In
  other words, we see that $G^{\op}\subset\Factor(\cG\actson Y)$. We show that
  this is actually an equality. Let $f:Y\rightarrow Y$ be a measure
  preserving quotient map that commutes with the action of $\cG$.
  \begin{step}
    Because $f$ commutes with all the $\psi_\Delta$, there is an
    injective map $\alpha:G\rightarrow G$ such that
    $f(x)_g=x_{\alpha(g)}$ for almost all $x\in Y$, and for all $g\in
    G$.
  \end{step}
  Fix a set $U\subset Y_0$ with $0<\nu_0(U)<1$, and consider the group
  $H=\{\Delta\in\Autmp(Y_0,\nu_0)\mid\Delta(U)=U\}$. Then we know that
  $H$ acts weakly mixingly on both $U$ and $Y_0\setminusb U$. In
  particular, for every finite set $I$, we know that the $H$-invariant
  subsets of $(Y_0,\nu_0)^I$ are precisely the disjoint unions of sets
  of the form $U^{I_1}\times (Y_0\setminusb U)^{I\setminusb I_1}$ for
  subsets $I_1\subset I$.

  Fix $i_0\in G$ and denote by $f_{i_0}:Y\rightarrow Y_0$ the
  composition of $f$ with the quotient map onto the $i_0$-component.
  Write $\widetilde U=f_{i_0}^{-1}(U)$.
  For a fixed finite subset $I\subset G$ and an element $y\in
  Y_0^{I}$, denote
  \begin{align*}
    \widetilde U_{I,y}&=\left\{\left.x\in Y_0^{G\setminusb I}\,\right\vert\, f(y,x)_{i_0}\in U\right\}\\
    \widetilde U_{I}&=\left\{\left.y\in Y_0^{I}\,\right\vert\, \nu_0^{G\setminusb
      I}(\widetilde U_{I,y})>0\right\}.
  \end{align*}
  Now it is clear that $\widetilde U\subset \widetilde U_I\times
  Y_0^{G\setminusb I}$. Moreover, $\widetilde U_I$ is a disjoint union of
  sets of the form $U^{I_1}\times (Y_0\setminusb U)^{I\setminusb I_1}$ for
  subsets $I_1\subset I$.

  We claim that there is an $i\in I$ such that $U^{\{i\}}\times
  Y_0^{I\setminus \{i\}}\subset \widetilde U_I$. To prove this claim, put
  $n=\nelt{I}$ and observe that the function $p\mapsto
  (1-p)^{n-1}+np(1-p)^{n-2}$ tends to $1$ as $p\rightarrow 0$, but the
  derivative in $0$ is positive. So, taking $p$ small enough, we can
  assume that $(1-p)^p+np(1-p)^{n-1}>(1-p)$. Take a subset $V\subset
  U$ with measure $\nu_0(V)=p$. Write $\widetilde
  V=f_{i_0}^{-1}(V)$ and define $\widetilde V_I$ in the same way as we
  defined $\widetilde U_I$. Then it is clear that
  $\widetilde V_I\subset \widetilde U_I$, but also that
  $\nu_0^I(\widetilde V_I)\geq p$. Observe that $\widetilde V_I$ is a
  disjoint union of sets of the form $V^{I_1}\times (Y_0\setminusb
  V)^{I\setminusb I_1}$ for some subset $I_1\subset I$. It follows
  that $\widetilde V_I$ must contain at least one of the sets
  $(Y_0\setminusb V)^I$ or $(Y_0\setminusb V)^{I\setminusb\{i\}}\times
  V^{\{i\}}$ for some $i\in I$. In the first case, it follows that
  $Y_0^I=\widetilde U_I$, while in the second case we find that
  $Y_0^{I\setminusb\{i\}}\times U^{\{i\}}\subset \widetilde U_I$. In both
  cases, we have proven our claim.

  Observe that, for any $\varepsilon>0$, there is a finite set
  $I\subset G$ such that $\nu(\widetilde U_{I}\times
  Y_0^{G\setminusb I}\setminusb \widetilde U)<\varepsilon$. Taking
  $\varepsilon = \nu_0(U)(1-\nu_0(U))$, we find a finite subset
  $I_0\subset G$ such that $\nu(\widetilde U_{I}\times
  Y_0^{G\setminusb I}\setminusb \widetilde U)<\varepsilon$. Then there is a unique $i\in I_0$ with
  $Y_0^{I_0\setminusb\{i\}}\times U^{\{i\}}\subset \widetilde
  U_{I_0}$. For every finite set $I_0\subset I\subset G$, we see that
  $Y_0^{I\setminusb\{i\}}\times U^{\{i\}}\subset \widetilde
  U_{I}$, with the same $i\in I_0\subset I$. As a consequence we find
  that $Y_0^{G\setminusb\{i\}}\times U^{\{i\}}\subset \widetilde
  U$. Comparing the measures, we see that this inclusion is actually
  an equality. In other words, $f_{i_0}^{-1}(U)=U^{\{i\}}$.

  Every subset $V\subset U$ is of the form $V=U\cap \Delta(U)$ for
  some automorphism $\Delta:Y_0\rightarrow Y_0$. Hence we compute that
  \begin{align*}
    f_{i_0}^{-1}(V)&=f^{-1}(U)\cap \psi_\Delta(f_{i_0}^{-1}(U))\\
    &=U^{\{i\}}\cap \Delta(U)^{\{i\}}\\
    &=V^{\{i\}}
  \end{align*}
  The same proof works for subsets of $Y_0\setminusb U$. Hence we get
  that for every $V\subset Y_0$, the inverse image $f_{i_0}^{-1}(V)$
  equals $V^{\{i\}}$. In other words, $f(x)_{i_0}=x_i$ almost everywhere.
  
  We can do this for every $i_0\in G$ and find a map
  $\alpha:G\rightarrow G$ such that $f(x)_g=x_{\alpha(g)}$. The map
  $\alpha$ is injective because otherwise $f$ would not be measure
  preserving. This finishes the proof of step 1.

  \begin{step}
    Because $f$ also commutes with the $\varphi_{U,g,\Delta}$, it
    follows that $f=f_k$ for some $k\in G$.
  \end{step}
  In fact, we only need that $f$ commutes with $\varphi_{U,g,\Delta}$
  for one non-trivial set $U\subset Y_0$ and one non-trivial
  automorphism $\Delta:Y_0\rightarrow Y_0$, but for all $g\in
  G$. Remember that $f$ is of the form $f(x)_h=x_{\alpha(h)}$ almost
  everywhere and for all $h\in G$. The fact that $f$ commutes with
  $\varphi_{U,g,\Delta}$ implies that $\alpha(hg)=\alpha(h)g$. Set
  $k=\alpha(e)$ where $e$ is the unit of $g$. Then it follows that
  $\alpha(g)=\alpha(e)g=kg$. In other words, we see that $f=f_k$.
  This concludes the proof of step 2, and hence the proof of lemma
  \ref{lem:discrete-factor}.
\end{proof}

Finally we prove the technical lemma we needed in the proof of theorem \ref{thm:ex-finite}.
\begin{lemma}
  \label{lem:no-fd}
  Let $\Lambda\actson (Y,\nu)$ be an action of a countable group on a probability space $(Y,\nu)$,
  then there is a probability measure preserving action $\widetilde\Lambda\actson (\widetilde Y,\tilde \nu)$
  of an anti-(T) group $\widetilde\Lambda$ (see definition \ref{def:anti-T}) such that 
  \begin{align*}
    \Factor(\widetilde\Lambda\actson(\widetilde Y,\tilde \nu))&=\Factor(\Lambda\actson(Y,\nu))
  \end{align*}
  and such that this new action does not have any non-trivial cocycles to compact groups. Moreover, the new action is ergodic.
\end{lemma}
\begin{proof}
  We use the generalized co-induced actions that were introduced in
  \cite{Deprez:Fundg}. For the convenience of the reader, we repeated the
  construction and basic properties in preliminary
  \ref{subsect:prelim:gencoind}.

  Without loss of generality, we can assume that $\Lambda=\FG_\infty$. Denote by $a_n$ the $n$-th canonical generator. Consider now the following group:
  \def\undereq#1_#2{\underset{\overset{\Vert}{#2}}{#1}}
  \begin{align*}
  \widetilde\Lambda&=\underbrace{\SL_2\IQ\ltimes\IQ^2}_H\free_{\Sigma}\underbrace{(\undereq{\SL_2\IZ}_{G_0}\ltimes(\undereq{\IZ^2}_A\times(\undereq{\IZ^2}_{B}\free
  \undereq{\IZ^2}_{F_1}\free\undereq{\IZ^2}_{F_2}\free\ldots)))}_{G},\\\text{ where }
  \Sigma&=\left\{\left(\begin{matrix}1&n\\0&1\end{matrix}\right)\mid n\in\IZ\right\}\subset G\\
  &\cong\left\{\left(\begin{matrix}1&-1\\1&0\end{matrix}\right)^n\mid n\in\IZ\right\}\subset H.
  \end{align*}
  Observe that $\widetilde \Lambda$ is an anti-(T) group because it is
  an amalgamated free product of poly-Haagerup groups.
  There is an obvious quotient from $\widetilde\Lambda$ onto the group
  \[\widetilde\Lambda_1=H\free_\Sigma(G_0\ltimes(A\times B)),\]
  where the $F_n$ are mapped to the identity element in
  $\widetilde\Lambda_1$. Consider the group
  $\widetilde\Lambda_0=\ST_2\IQ$ of upper triangular matrices in
  $\SL_2\IQ$, and let $\widetilde\Lambda$ act on
  $I=\widetilde\Lambda_1/\widetilde\Lambda_0$
  by left translation.
  Define a cocyle $\omega:\widetilde\Lambda\times I\rightarrow \FG_\infty$ by the following relations.
  \begin{align*}
    \omega(s,i)&=e&&\text{ if }s\in\widetilde\Lambda_0\\
    \omega(s,i)&=e&&\text{ if }s\in F_n\text{ and a reduced word for }i\text{ starts with a letter from }H\\
    \omega(s,i)&=a_n^{\det(s\vert b)}&&
    \begin{aligned}&\text{ if }s\in F_n\text{ and
    }\\&(a,b)g\in G_0\ltimes(A\times B)\text{ is the first letter of a reduced word for }i.
    \end{aligned}
  \end{align*}
  Above, we denoted $\det(s\vert b)$ for the determinant of the matrix whose columns are $s$ and $b$.

  Consider the generalized co-induced action
  $\widetilde\Lambda\actson\widetilde Y$ of $\Lambda\actson Y$,
  associated to the cocycle $\omega$. This cocycle clearly satisfies the conditions of lemma
  \ref{lem:gen-co-ind}. So we see already that
  \[\Factor(\widetilde\Lambda\actson(\widetilde Y,\tilde
  \nu))=\Factor(\Lambda\actson(Y,\nu)),\]
  and moreover that $\widetilde\Lambda$ acts ergodically on
  $\widetilde Y$.

  It remains to show that every cocycle
  $\alpha:\widetilde\Lambda\times\widetilde Y\rightarrow \cG$, to a
  compact group $\cG$, is in fact trivial. When restricted to
  $\widetilde\Lambda_0$, the action is just a generalized Bernoulli
  action. We can apply Popa's cocycle superrigidity theorem
  \cite{Popa:CocycleSuperrigidityMalleable} to the relatively rigid inclusion $\IZ^2\subset\IQ^2\subset H$.
  We find a measurable map $\varphi:\widetilde
  Y\rightarrow \cG$ such that
  $\varphi(gx)\alpha(g,x)\varphi(x)^{-1}$ is independent of the
  $x$-variable for every $g\in\IZ^2$. Since $\IZ^2$ is almost normal
  in $H$ and acts weakly mixingly on $\widetilde Y$, the same is in
  fact true for all $g\in H$.

  Analogously, using the rigid inclusion $A\subset G$, we find a measurable function
  $\varphi_2:\widetilde Y\rightarrow \cG$ such that
  $\varphi_2(gx)\alpha(g,x)\varphi_2(x)^{-1}$ is independent of the
  $x$-variable for all $g\in G$. For $g\in\Sigma$, both
  \[\theta(g)=\varphi(gx)\alpha(g,x)\varphi(x)^{-1}\text{ and
  }\theta_2(g)=\varphi_2(gx)\alpha(g,x)\varphi_2(x)^{-1}\]
  are independent of $x$. So we find that
  \[(\varphi\varphi_2^{-1})(gx)=\theta(g)(\varphi\varphi_2^{-1})(x)\theta_2(g)^{-1}\]
  for all $g\in \Sigma$. Since $\Sigma$ acts weakly mixingly
  on $\widetilde Y$, it follows that $\varphi=\varphi_2$ up to a
  constant. So we can assume that $\varphi$ actually equals
  $\varphi_2$, and $\alpha$ is cohomologuous to a group morphism
  $\theta:\widetilde\Lambda\rightarrow \cG$.

  Since the only group morphism $H\rightarrow \cG$ is the trivial
  morphism, we see that at least $H\subset\ker \theta$. In particular,
  $\Sigma\subset \ker\theta$. The smallest normal subgroup of $G$ that
  contains $\Sigma$ is $G$ itself, so it follows that $G\subset
  \ker\theta$. But then $\theta$ is the trivial group morphism.
\end{proof}


\newcommand{\etalchar}[1]{$^{#1}$}
\bibliography{references}{}
\bibliographystyle{sdpabbrv}
\end{document}

%% file: sdpDefinitions.tex
\usepackage{amsmath}
\usepackage{amssymb}
\usepackage{amsfonts}
\usepackage{amsthm}
\usepackage{amscd}
\usepackage{tikz}
\usepgflibrary{arrows}
\newtheorem{theorem}{Theorem}[section]
\newtheorem{lemma}[theorem]{Lemma}

\newtheorem{definition}[theorem]{Definition}

\newtheorem{observation}[theorem]{Observation}
\newcounter{step}[theorem]
\newenvironment{step}{\par\refstepcounter{step}\textbf{step \thestep: }\begingroup\it}{\endgroup\par}
\newcommand\len[1]{\left\lvert #1\right\rvert}
\newcommand\norm[1]{\left\lVert #1\right\rVert}
\newcommand\nelt[1]{\left\lvert #1\right\rvert}
\newcommand\abs[1]{\left\lvert #1\right\rvert}

\DeclareMathOperator{\IC}{\mathbb{C}}
\DeclareMathOperator{\IR}{\mathbb{R}}
\DeclareMathOperator{\IRpos}{\mathbb{R}_{+}^{\ast}}
\DeclareMathOperator{\IQ}{\mathbb{Q}}
\DeclareMathOperator{\IZ}{\mathbb{Z}}
\DeclareMathOperator{\IN}{\mathbb{N}}
\DeclareMathOperator{\IH}{\mathbb{H}}
\DeclareMathOperator{\F}{F}

\DeclareMathOperator{\cF}{\mathcal{F}}

\DeclareMathOperator{\cU}{\mathcal{U}}
\DeclareMathOperator{\cG}{\mathcal{G}}

\def\MatrixGroup#1_#2#3{\mathchoice{#1_{#2}\hspace*{-0.3mm}#3}{#1_{#2}\hspace*{-0.4mm}#3}{#1_{#2}\hspace*{-0.4mm}#3\hspace*{0.3mm}}{#1_{#2}\hspace*{-0.4mm}#3\hspace*{0.3mm}}%
  \def\tempa{#3}%
  \def\tempb{\IQ}\ifx\tempa\tempb\else%
  \def\tempb{\IZ}\ifx\tempa\tempb\else%
  \def\tempb{R}\ifx\tempa\tempb\else%
  \def\tempb{(\F_2[X])}\ifx\tempa\tempb\else%
  \def\tempb{(\F_2)}\ifx\tempa\tempb\else%
  \def\tempb{\IC}\ifx\tempa\tempb\else%
  \fout\fi\fi\fi\fi\fi\fi%
}

\def\ST{\MatrixGroup\OpST}
\def\SL{\MatrixGroup\OpSL}
\def\GL{\MatrixGroup\OpGL}
\DeclareMathOperator{\OpGL}{GL}

\DeclareMathOperator{\OpST}{ST}
\DeclareMathOperator{\OpSL}{SL}
\DeclareMathOperator{\Sp}{Sp}
\DeclareMathOperator{\liesl}{\mathfrak{sl}}
\DeclareMathOperator{\liesp}{\mathfrak{sp}}
\DeclareMathOperator{\PSp}{PSp}
\DeclareMathOperator{\Hurwitz}{\mathbb{H}ur}
\newcommand{\FG}{\mathbb{F}}
\DeclareMathOperator{\free}{\ast}

\makeatletter
\def\quot{\@ifnextchar[{\sdp@quotleftA}{\sdp@quotright}}
\def\sdp@quotleftA[#1]#2{\@ifnextchar[{\sdp@quotboth[{#1}]{#2}}{\sdp@quotleft[{#1}]{#2}}}
\def\sdp@quotright#1[#2]{#1/#2}
\def\sdp@quotboth[#1]#2[#3]{#1\backslash #2 / #3}
\def\sdp@quotleft[#1]#2{#1\backslash #2}
\makeatother

\DeclareMathOperator{\Centr}{Centr}

\DeclareMathOperator{\Comm}{Comm}

\DeclareMathOperator{\Stab}{Stab}

\DeclareMathOperator{\RFBimod}{RFBimod}

\DeclareMathOperator{\Aut}{Aut}
\DeclareMathOperator{\End}{End}
\DeclareMathOperator{\Factor}{Factor}

\def\Autmp(#1,#2){\Aut_{#2}(#1)}
\def\Auttp(#1,#2){\Aut_{#2}(#1)}

\DeclareMathOperator{\Out}{Out}
\def\Outmp(#1,#2){\Out_{#2}(#1)}
\def\Outtp(#1,#2){\Out_{#2}(#1)}
\DeclareMathOperator{\Unitary}{\mathcal{U}}
\DeclareMathOperator{\Fundg}{\mathcal{F}}

\DeclareMathOperator{\op}{op}

\def\Rel(#1\actson #2){\mathcal{R}(#1\actson #2)}

\def\fullg(#1){\left[#1\right]}
\def\fullpg(#1){\left[\left[#1\right]\right]}

\DeclareMathOperator{\Centre}{\mathcal{Z}}
\DeclareMathOperator{\vnNorm}{\mathcal{N}}
\DeclareMathOperator{\Bounded}{\mathcal{B}}
\DeclareMathOperator{\MatM}{M}
\DeclareMathOperator{\MatD}{D}
\DeclareMathOperator{\MatU}{U}
\DeclareMathOperator{\Lp}{L}
\DeclareMathOperator{\Lg}{L}

\DeclareMathOperator{\E}{E}
\DeclareMathOperator{\Tr}{Tr}

\newcommand{\bimod}[3]{\strut_{#1}#2\strut_{#3}}
\DeclareMathOperator*{\embeds}{\prec}
\DeclareMathOperator*{\nembeds}{\not\prec}
\DeclareMathOperator*{\fembeds}{\prec^f}

\def\vnInterB(#1){\zH_{#1}}
\DeclareMathOperator{\Ad}{Ad}
\DeclareMathOperator{\actson}{\curvearrowright}
\DeclareMathOperator{\id}{id}

\DeclareMathOperator{\setminusb}{--\,}

\newcommand\restrict[1]{\vert_{#1}}
\makeatletter
\def\zH{\@ifnextchar_{\zH@}{\@undefined}}
\DeclareMathOperator{\zH@}{H}
\def\smallscripts#1{\@ifnextchar_{\smallscripts@{#1}}{\@undefined}}
\def\smallscripts@#1{\@ifnextchar_{\smallscripts@subP{#1}}{\@ifnextchar^{\smallscripts@superP{#1}}{\smallscripts@none{#1}}}}
\def\smallscripts@subP#1_#2{\@ifnextchar^{\smallscripts@subsuperP{#1}{#2}}{\smallscripts@sub{#1}{#2}}}
\def\smallscripts@subsuperP#1#2^#3{\smallscripts@both{#1}{#2}{#3}}
\def\smallscripts@superP#1^#2{\@ifnextchar_{\smallscripts@supersubP{#1}{#2}}{\smallscripts@super{#1}{#2}}}
\def\smallscripts@supersubP#1#2_#3{\smallscripts@both{#1}{#3}{#2}}
\def\smallscripts@none#1{{\displaystyle #1} }
\def\smallscripts@sub#1#2{{\displaystyle #1}_{\scriptscriptstyle #2}}
\def\smallscripts@super#1#2{{\displaystyle #1}^{\scriptscriptstyle #2}}
\def\smallscripts@both#1#2#3{{\displaystyle #1}_{\scriptscriptstyle #2}^{\scriptscriptstyle #3}}
\makeatother
\newcommand\chara{\smallscripts{\chi}}
